%% file: root.tex
\documentclass[twocolumn]{IEEEtran}
\usepackage{cite}
\usepackage{color,yfonts}
\usepackage{graphicx}
\graphicspath{{fig/}{jpeg/}}
\usepackage[cmex10]{amsmath}
\usepackage{amssymb}
\usepackage{epstopdf}
\usepackage{multirow}
\usepackage{algorithm}
\usepackage{tikz}
\usetikzlibrary{shapes,arrows,matrix}
\usetikzlibrary{decorations.pathmorphing} 
\usetikzlibrary{fit}					
\usetikzlibrary{backgrounds}	
\usepackage{verbatim}
\usepackage{algpseudocode}
\usepackage{amsmath,amssymb,lipsum}
\usepackage{hyperref}
\usepackage{comment}
\usepackage{graphicx}
\usepackage[font=small]{caption}
\usepackage{subcaption}
\usepackage{mathtools,enumerate}
\input{mysymbol.sty}
\input{my_sections.tex}

\usepackage{theorem}
\newtheorem{thm}{Theorem}
\newtheorem{lemma}{Lemma}

\theoremstyle{definition}
\newtheorem{assumption}{Assumption}
\newtheorem{remark}{Remark}

\def \QED {$\blacksquare$}

\def \x {{\mathbf{x}}}

\def \t {{\mathbf{t}}}
\def \u {{\mathbf{u}}}
\def \y {{\mathbf{y}}}
\def \h {{\mathbf{h}}}
\def \ep {{\epsilon}}

\newcommand{\qed}{\hfill\blacksquare}
 \providecommand{\Ex}[1]{\mathbb{E}\left[#1\right]}
 \providecommand{\abs}[1]{\left|#1\right|}
 \providecommand{\norm}[1]{\left\|#1\right\|}
 \providecommand{\ip}[1]{\boldsymbol{\langle}#1\boldsymbol{\rangle}}
 
  \newcommand{\col}[1]{\textcolor{black}{#1}}
  \newcommand{\colb}[1]{\textcolor{blue}{#1}}
    
  \tikzstyle{agent}=[circle,
  thick,
  minimum size=.5cm,
  draw=blue!80,
  fill=blue!20]

  \tikzstyle{neighbor}=[circle,
  thick,
  minimum size=.5cm,
  draw=cyan!85!black,
  fill=cyan!40,
  decorate,
  decoration={random steps,
  	segment length=2pt,
  	amplitude=2pt}]
  
  \tikzstyle{local_nat}=[rectangle,
  thick,
  minimum size=.5cm,
  draw=red!80,
  fill=red!20]
  
  \tikzstyle{glob_nat}=[rectangle,
  thick,
  minimum size=.5cm,
  draw=red!100,
  fill=red!40]      
  
  \tikzstyle{background}=[rectangle,
  fill=gray!10,
  inner sep=0.2cm,
  rounded corners=5mm]
  
  \tikzstyle{background2}=[rectangle,
  fill=green!20,
  inner sep=0.2cm,
  rounded corners=5mm]

\title{\vspace{-1cm}\blue{Asynchronous Decentralized Stochastic Optimization in Heterogeneous Networks}}
\author{Amrit~Singh~Bedi$^\star$,~\IEEEmembership{Student Member,~IEEE,}
	Alec~Koppel$^\dagger$,~\IEEEmembership{Student Member,~IEEE,}
	and~Ketan~Rajawat$^\star$,~\IEEEmembership{Member,~IEEE}
	\thanks{
		A. S. Bedi and K. Rajawat are with the Department of Electrical Engineering,
		Indian Institute of Technology Kanpur, Kanpur 208016, India (e-mail:
		amritbd@iitk.ac.in; ketan@iitk.ac.in). A. Koppel is with the U.S. Army Research Laboratory, Adelphi, MD, USA. (e-mail: akoppel@seas.upenn.edu.). }
	\vspace{-1cm}	}
\begin{document}
\maketitle

\begin{abstract}
	We consider expected risk minimization in multi-agent systems comprised of distinct subsets of agents operating without a common time-scale. 
	Each individual in the network is charged with minimizing the global objective function, which is an average of sum of the statistical average loss function of each agent in the network. Since agents are not assumed to observe data from identical distributions, the hypothesis that all agents seek a common action is violated, and thus the hypothesis upon which consensus constraints are formulated is violated. Thus, we consider nonlinear network proximity constraints which incentivize nearby nodes to make decisions which are close to one another but not necessarily coincide. Moreover, agents are not assumed to receive their sequentially arriving observations on a common time index, and thus seek to learn in an asynchronous manner. An asynchronous stochastic variant of the Arrow-Hurwicz saddle point method is proposed to solve this problem which operates by alternating primal stochastic descent steps and Lagrange multiplier updates which penalize the discrepancies between agents. This tool leads to an implementation that allows for each agent to operate asynchronously with local information only and message passing with neighbors. Our main result establishes that the proposed method yields convergence in expectation both in terms of the primal sub-optimality and constraint violation to radii of sizes $\ccalO(\sqrt{T})$ and $\ccalO(T^{3/4})$, respectively. Empirical evaluation on an asynchronously operating wireless network that manages user channel interference through an adaptive communications pricing mechanism demonstrates that our theoretical results translates well to practice. 
\end{abstract}

\section{Introduction}\label{sec:intro}
In emerging technologies such as wireless communications and networks consisting of interconnected consumer devices \cite{App_3}, increased sensing capabilities are leading to new theoretical challenges to classical parameter estimation. These challenges include the fact that data is persistently arriving in a sequential fashion \cite{koppel2016online}, that it is physically decentralized across an interconnected network, and that the nodes of the network may correspond to disparate classes of objects (such as users and a base station) with different time-scale requirements \cite{resource_alloca_tut_4}. In this work, we seek to address this class of problems through extensions of online decentralized convex optimization \cite{Shalev-Shwartz2012} to the case where the agents of the network may be of multiple different classes, and operate on  different time-scales \cite{Tsitsiklis1986}.

To address the fact that we seek iterative tools for streaming data, we consider stochastic optimization problems \cite{Vapnik1995,bousquet2008tradeoffs}. In this setting, the objective function $\mathbb{E}[f(\bbx,\bbtheta)]$ is an expectation over a set of functions parameterized by a random variable $\bbtheta$. The objective function encodes, for example, the quality of a statistical parameter estimate. Through a sequence of realizations of a random variable $\bbtheta_t$, we seek to find parameters that are good with respect to the average objective. The classical method to address this problem is stochastic gradient descent (SGD), which involves descending along the negative of the stochastic gradient in lieu of the true gradient to circumvent the computation of an infinite complexity expectation \cite{Robbins1951,Bottou1998}. SGD forms the foundation of tools considered in this paper for asynchronous multi-agent settings.


Here we seek solutions to stochastic programs in which data is scattered across an interconnected network $\ccalG\!\!=\!\!(\!\ccalV,\! \ccalE\!)$ of agents, each of which is associated with a unique stream of data $\{\bbtheta^i_t\}_{t\geq0}$. Agents $i\!\in\!\ccalV$ then seek to find a solution based on local computations only which is as good as one at a centralized location: this setting is mathematically defined by introducing a local copy of a global parameter estimate and then having each agent seek to minimize a global sum of all local objectives $\sum_i \mathbb{E}[f^i(\bbx^i,\bbtheta^i)]$ while satisfying consensus constraints  $\bbx_i = \bbx_j$ for all node pairs $(i,j) \in \ccalE$. 
Techniques that are good for distributed convex optimization, for example, those based on penalty method  \cite{Nedic2009a,tsianos2012distributed} or on Lagrange duality \cite{nedic2009subgradient,boyd2011distributed}, have in most cases translated into the distributed stochastic domain without major hurdles, as in \cite{tsianos2012distributed,Koppel2015a,rabbat2014asynchronous,ling2014decentralized}.

{In the aforementioned works, consensus constraints are enforced in order to estimate a common decision variable while leveraging parallel processing architectures to achieve computational speedup \cite{bertsekas1989parallel,recht2011hogwild}. }
Contrariwise, when unique priors on information available at distinct group of agents are available as in sensor \cite{koppel2017proximity} or robotic \cite{koppel2017d4l} networks, enforcing consensus degrades the statistical accuracy of each agent's estimate \cite{koppel2017proximity}. Specifically, if the observations at each node are independent but \emph{not} identically distributed, consensus may yield a sub-optimal solution. Motivated by heterogeneous networked settings \cite{koppel2017d4l,bedi2017asynchronous} where each node observes a unique local data stream, we focus on the setting of multi-agent stochastic optimization with nonlinear \emph{proximity} constraints which incentivize nearby nodes to select estimates which are similar but not necessarily equal.

In the setting of nonlinear constraints, penalty methods such as distributed gradient descent do not apply \cite{Nedic2009a}, and dual or proximal methods require a nonlinear minimization in an inner-loop of the algorithm \cite{boyd2011distributed}. Therefore, we adopt a method which hinges on Lagrange duality that avoids costly $\argmin$ computations in the algorithm inner-loop, namely primal-dual method \cite{Arrow1958}, also referred to as saddle point method. Alternative attempts to extend multi-agent optimization techniques to heterogeneously correlated problems have been considered in \cite{Caruana:1997:ML:262868.262872,chen2014multitask} for special loss functions and correlation models, but the generic problem was online recently solved in \cite{koppel2017proximity} with a stochastic variant of the saddle point method. 

However, insisting on all agents to operate on a common clock creates a bottleneck for implementation in practical settings because typically nodes may be equipped with different computational capacity due to power and energy design specifications, as well as a difference in the sparsity of each agent's data stream.  Therefore, we attempt to extend multi-agent stochastic optimization with nonlinear constraints to asynchronous settings \cite{duchi2015asynchronous}. Asynchrony in online optimization has taken on different forms, such as, for instance, maintaining a local Poisson clock for each agent \cite{srivastava2011distributed } or a distribution-free generic bounded delay \cite{Tsitsiklis1986,sayed2015asynchronous}, the approach considered here.

In this paper, we extend the primal-dual method of \cite{Arrow1958,koppel2017proximity} for multi-agent stochastic optimization problems with nonlinear network proximity constraints to asynchronous settings. The proposed algorithm allows the gradient to be delayed for the primal and dual updates of the saddle point method. 
The main technical contribution of this paper is to provide mean convergent results for both the global primal cost and constraint violation, establishing that the Lyapunov stability results of \cite{koppel2017proximity} translate successfully to asynchronous computing architectures increasingly important in intelligent communication systems. Empirical evaluation on an asynchronously operating cellular network that manages cross-tier interference through an adaptive pricing mechanism demonstrates that our theoretical results translates well to
practice.

\colb{The rest of the paper is organized as follows. A multi-agent optimization problem without consensus is formulated in Section  \ref{sec:prob}. An asynchronous saddle point algorithm is proposed to solve the problem in Section \ref{sec:algorithm}.  The detailed convergence analysis for the proposed algorithm is presented in Section  \ref{sec:convergence}. Next, a practical problem of interference management through pricing is solved in Section \ref{sec:fields}. Section \ref{sec:conclusion} concludes the paper. }
\section{ Multi-Agent Optimization without Consensus}\label{sec:prob}

We consider agents $i$ of a symmetric, connected, and directed network $\ccalG = (V, \ccalE)$ with $|V|=N$ nodes and $|\ccalE|=M$ edges. Each agent is associated with a (non-strongly) convex loss function $f^i : \ccalX \times \Theta_i \rightarrow \reals$ that is parameterized by a $p$-dimensional decision variable $\bbx^i \in \ccalX\subset \reals^p$ and a random vector $\bbtheta_i \in \Theta_i\subset \reals^q$. The functions $f^i(\bbx^i,{\bbtheta^i})$ for different ${\bbtheta^i}$ encodes the merit of a particular linear statistical model $\bbx^i$, for instance, and the random vector $\bbtheta$ may be particularized to a random pair $\bbtheta=(\bbz, \bby)$. In this setting, the random pair corresponds to feature vectors $\bbz$ together with their binary labels $\bby \in \{-1, 1\}$ or real values $\bby \in \reals$, for the respective problems of classification or regression.
Here we address the case that the local random vector ${\bbtheta^i}$ represents data which is revealed to node $i$ \emph{sequentially} as realizations $\bbtheta^{i}_t$ at time $t$, and agents would like to process this information on the fly. Mathematically this is equivalent to the case where the total number of samples $T$ revealed to agent $i$ is not necessarily finite. A possible goal for agent $i$ is the solution of the local expected risk minimization problem, 
\begin{equation} \label{eq:dist_stoch_opt_local}
\bbx^{\text{L}}(i) := \argmin_{\bbx^i \in \ccalX} 
F^i(\bbx^i)
:= \argmin_{\bbx^i\in \reals^p} 
\mbE_{{\bbtheta^i}}[f^i(\bbx^i,{\bbtheta^i})]\; . 
\end{equation}
where we define $F^i(\bbx^i):=\mbE_{{\bbtheta^i}}[f^i(\bbx^i,{\bbtheta^i})]$ as the local average function at node $i$. We also restrict $\ccalX$ to be a compact convex subset of $\reals^p$ associated with the $p$-dimensional parameter vector of agent $i$. By stacking the problem \eqref{eq:dist_stoch_opt_local} across the entire network, we obtain the equivalent problem
\begin{equation} \label{eq:dist_stoch_opt}
\bbx^{\text{L}} = \argmin_{\bbx\in \ccalX^N} F(\bbx) 
:= \argmin_{\bbx \in \ccalX^N} \sum_{i=1}^N \mbE_{{\bbtheta^i}}[f^i(\bbx^i,{\bbtheta^i})]\; . 
\end{equation}
where we define the stacked vector $\bbx=[{\bbx^1,\ldots,\bbx^N}]\in\ccalX^N\subset \reals^{Np}$, and the global cost function $F(\bbx) := \sum_{i=1}^N \mbE_{{\bbtheta^i}}[f^i(\bbx^i,{\bbtheta^i})]$. We define the global instantaneous cost similarly:  $f(\bbx, \bbtheta)=\sum_i f^i(\bbx^i, {\bbtheta^i})$. 

Note that \eqref{eq:dist_stoch_opt_local} and \eqref{eq:dist_stoch_opt} describe the same problem since the variables $\bbx^i$ at different agents are not coupled to one another. In many situations, the parameter vectors of distinct agents are related, and thus there is motivation to couple the estimates of distinct agents to each other such that one agent may take advantage of another's data. Most distributed optimization works, for instance, consensus optimization, hypothesize that all agents seek to learn the common parameters $\bbx^i$ for all $i\in V$, i.e.,
$    {\bbx^i = \bbx^j}, \text{ for all } j\in n_i \; .$
where $n_i$ denotes the neighborhood of agent $i$. Making all agents variables equal only makes sense when agents observe information drawn from a common distribution, which is the case for industrial-scale machine learning, but is predominantly not the case for sensor \cite{koppel2017proximity} and robotic networks \cite{koppel2017d4l}. As noted in \cite{koppel2017proximity}, generally, nearby nodes observe similar but not identical information, and thus to incentivize collaboration without enforcing consensus, we introduce a convex local proximity function with real-valued range of the form $h^{ij}(\bbx^i, \bbx^{j},\bbtheta^i,\bbtheta^j )$  that depends on the observations of neighboring agents and a tolerance $\gamma_{ij}\geq 0$. These stochastic constraints then couple the decisions of agent $i$ to those of its neighbors $j \in {n_i}$ as the solution of the constrained stochastic program
\begin{alignat}{2} \label{eq:coop_stoch_opt}
\bbx^*\in \ 
&\argmin_{\bbx \in \ccalX^{N} }\ &&\sum_{i=1}^N \mbE_{\bbtheta_i}[f^i(\bbx_i,\bbtheta_i)] \\
&\text{\ s.t. } \ &&\!\!\!\!\mathbb{E}_{\bbtheta^i, \bbtheta^j}\left[h^{ij}({\bbx^i, \bbx^{j},\bbtheta^i,\bbtheta^j })\right]              \leq \gamma_{ij}, \text{ for all }j\in n_i.\nonumber
\end{alignat}
%
Examples of the constraint in the above formulation include approximate consensus constraints $\norm{\x^i-\x^j} \leq \gamma_{ij}$, quality of service  $\textbf{SINR}(\x^i,\x^j) \geq \gamma_{ij}$ where SINR is the signal-to-interference-plus-noise function, relative entropy $D(\x^i\mid\mid \x^j) \leq \gamma_{ij}$, or   budget $\gamma_{ij}^{\min} \leq x^i + x^j \leq \gamma_{ij}^{\max}$ constraints. In this work, we seek decentralized online solutions to the constrained problem \eqref{eq:coop_stoch_opt} \emph{without} the assumption that agents operate on a common time index, motivated by the fact that asynchronous computing settings are common in large distributed wireless networks. 
In the next section we turn to developing an algorithmic solution that meets these criteria.

%
\begin{remark}\label{remark1}
	\normalfont  The pairwise stochastic constraints in \eqref{eq:coop_stoch_opt} can be readily generalized to arbitrary neighborhood constraints:
	\vspace{-2mm}
	\begin{align}\label{eq:const_gen}
	\Ex{\h^i(\{\x^j, \bbtheta^j\}_{j \in n_i'})} \leq \mathbf{0}
	\end{align}
	where the set $n_i':=n_i \cup \{i\}$ denotes the set of all neighbors of node $i$ inlcuding the node $i$ itself. It can be seen that constraint in \eqref{eq:coop_stoch_opt} is a special case of that in \eqref{eq:const_gen}, with the $j$-th entry of the $\abs{n_i} \times 1$ vector function $\h{^i}(\cdot)$ defined as
	\begin{align}\label{constraint_gen}
	\left[\h{^i}(\{\x^j, \bbtheta^j\}_{j \in n_i'})\right]_j := h^{ij}(\x^i,\x^j,\bbtheta^i,\bbtheta^j) - \gamma_{ij}.
	\end{align}
	More generally, \eqref{eq:const_gen} allows the consensus constraints to be imposed on the entire neighborhood of a node. For instance the approximate version of the consensus constraint $\x^i = (1/{|n_i|} )\sum_{j \in n_i} \x^j$ takes the form \vspace{-2mm}
	\begin{align}
	\big\|\x^i - (1/{|n_i|}) \sum_{j \in n_i} \x^j\big\| \leq \gamma_{i}.
	\end{align}
	Such general constraints also arise in communication systems in form of SINR constraints. For instance, consider a communication system where the interference at node $i$ from $j$ can be written as $p^j(\x^j,g^{ij})$ with $g^{ij}$ denoting the channel gain from node $j$ to node $i$. Then the SINR constraint at node $i$ is of the form in \eqref{eq:const_gen}, and is given by 
	\begin{align}
	h^i(\{\x^j, \bbtheta^j\}_{j \in n_i'}) = \gamma_{ij} - \frac{p^i(\x^i,g^{ii})}{\sigma^2 + \sum_{j\in n_i} p^j(\x^j,g^{ij})} 
	\end{align}
	where $\sigma^2$ denotes the noise power. Hence, the generalized stochastic problem can be expressed as 
	\begin{alignat}{3} \label{eq:coop_stoch_opt_gen}
	\bbx^*\in \ 
	&\argmin_{\bbx \in \ccalX^{N} }\ &&\sum_{i=1}^N \mbE[f^i(\bbx_i,\bbtheta_i)] \\
	&\text{\ s.t. } \ &&\!\!\!\!\Ex{\h^i\left({\{\bbx^{j},\bbtheta^j\}_{j\in n_i'}}\right)}\leq \boldsymbol{0}\ \ \text{for all $i$}.\nonumber
	\end{alignat}
	It is mentioned that the convergence analysis is performed for the generalized problem in \eqref{eq:coop_stoch_opt_gen} for clarity of exposition since in this more general setting we may vectorize the constraints, while the main results in Section \ref{sec:convergence} are presented for simpler problem of \eqref{eq:coop_stoch_opt} [see Appendix 0] for increased interpretability.  

\end{remark}
\section{Asynchronous Saddle Point Method}\label{sec:algorithm}
Methods based upon distributed gradient descent and penalty methods more generally \cite{DBLP:journals/corr/abs-1112-2972,RamNedicVeeravalli, Yuan2013} are inapplicable to settings with nonlinear constraints, with the exception of \cite{towfic2015stability}, which requires attenuating learning rates to attain constraint satisfaction. On the other hand, the dual methods proposed in \cite{1506308,Jakubiec2013,Wei2012} require a nonlinear minimization computation at each algorithm iteration, and thus is impractically costly. Therefore, in this section we develop a computationally light weight method based on primal-dual method that may operate in decentralized online asynchronous settings with constant learning rates that are better suited to changing environments.

For a decentralized algorithm, each node $i$ can access the information from its neighbors  $j\in n_i$ only. For an online algorithm, the stochastic i.i.d quantities of unknown distribution are observed sequentially $\bbtheta_t^i$  at each time instant $t$. In addition to these properties, an algorithm is called asynchronous, if parameter updates may be executed with out-of-date information and the requirement that distinct nodes operate on a common time-scale is omitted. To develop an algorithm which meets these specifications, begin by considering the approximate Lagrangian relaxation of \eqref{eq:coop_stoch_opt} stated as
\begin{align} \label{eq:lagrangian}
\ccalL(\bbx,\bblambda)= &\sum_{i=1}^N \Bigg[
\mbE\bigg[f^i(\bbx^i,\bbtheta^i) \\
&+
\sum_{j\in n_i}\lam^{ij} \left(h^{ij}\left(\bbx^{i},\bbx^{j},\bbtheta^i,\bbtheta^j \right)-\gamma_{ij}\right)
-\frac{\delta \eps}{2} (\lam^{ij})^2\bigg] \Bigg],\nonumber
\end{align}
where $\lam^{ij}$ is a non-negative Lagrange multiplier associated with the non-linear constraint in \eqref{eq:coop_stoch_opt}. Here,  $\bblam$ defines the collection of all dual variables $\lam^{ij}$ into a single vector $\bblam$. Observe that \eqref{eq:lagrangian} is not the standard Lagrangian of the \eqref{eq:coop_stoch_opt} but instead an augmented Lagrangian due to the presence of the term $-(\delta \eps/2)(\lam^{ij})^2$. This terms acts like a regularizer on the dual variable with associated parameters $\delta$ and $\eps$ that allow us to control the accumulation of constraint violation of the algorithm over time, as is discussed in the following section and proofs in the appendices.

The stochastic saddle point algorithm, when applied to \eqref{eq:lagrangian}, operates by alternating primal and dual stochastic gradient descent and ascent steps, respectively. We consider the stochastic saddle point method as a template upon which we construct an asynchronous protocol. Begin then by defining the stochastic approximation of the augmented Lagrangian evaluated at observed realizations $\bbtheta^{i}_t$ of the random vectors $\bbtheta^{i}$ for each $i\in\ccalV$:
\begin{align} \label{eq:stoch_lagrangian}
\hat{\ccalL}_t(\bbx,\bblambda)=& \sum_{i=1}^N \Big[f^i(\bbx^i,\bbtheta^i_t) \\
&+
\sum_{j\in n_i}\!\!\lam^{ij}\!\! \left(h^{ij}\left(\bbx^{i},\bbx^{j},\bbtheta^i_t,\bbtheta^j_t \right)\!-\!\gamma_{ij}\right)
\!-\!\frac{\delta \eps}{2} (\lam^{ij})^2\Big].\nonumber
\end{align}
The stochastic saddle point method applied to the stochastic Lagrangian \eqref{eq:stoch_lagrangian} takes the following form similar to \cite{koppel2017proximity} as
\vspace{-2mm}
\begin{align} \label{eq:sp_primal}
\bbx_{t+1}   &=\ccalP_{\ccalX}\Big[ \bbx_t - \eps \nabla_\bbx \hat{\ccalL}_t (\bbx_t, \bblam_t)\Big] \;,  \\
\bblam_{t+1} &= \Big[\bblam_t + \eps \nabla_{\bblam} \hat{\ccalL}_t (\bbx_{t}, \bblam_t)  \Big]_+\; ,
\label{eq:sp_dual}
\end{align}
where  $\col{\nabla_\bbx \hat{\ccalL} (\bbx_t, \bblam_t)}$ and $\nabla_{\bblam} \hat{\ccalL} (\bbx_t, \bblam_t)$, are the primal and dual stochastic  gradients\footnote{Note that these may be subgradients if the objective/ constraint functions are non-differentiable. The proof is extendable  to non-differentiable cases.  } of the augmented Lagrangian with respect to $\bbx$ and $\bblam$, respectively. These are not the actual gradients of \eqref{eq:lagrangian} rather are stochastic gradients calculated at the current realization of the random vectors $\bbtheta_t^i$ for all $i$. The component wise projection for a vector $\x$ on to the given compact set $\mathcal{X}$ is here denoted by $\ccalP_\ccalX(\bbx)$. Similarity,  $[ \cdot ]_+$ represents the component wise projection on to the positive orthant $\reals^M_{+}$. An important point here is that the method stated in \eqref{eq:sp_primal} - \eqref{eq:sp_dual} can be implemented with decentralized computations across the network, as stated in \cite[Proposition 1]{koppel2017proximity}. Here $\eps>0$ is a constant positive step-size.

Observe that the implementation of the \cite[Algorithm 1]{koppel2017proximity}, which is defined by \eqref{eq:sp_primal}-\eqref{eq:sp_dual}, it is mandatory to perform the primal and dual updates at each node with a common time index $t$. The update of primal variable at node $i$ requires the current gradient of its local objective function $\nabla_{\bbx^i} f^i (\bbx^{i}_t , \bbtheta^{i}_t )$ and current gradient from all the neighbors $j\in n_i$ of node $i$ as $\nabla_{\bbx^i} h^{ij}  \left(\bbx^{i}_t,\bbx^{j}_t,\bbtheta^{i}_t,\bbtheta^{j}_t\right)$. This availability of the gradients from the neighbors on a common time-scale is a strong assumption that insists upon perfect communications, similarity of computational capability of distinct nodes, and similar levels of sparsity among agents' data that are oftentimes violated in large heterogeneous systems. This limitation of synchronized methods motivates the subsequent development of asynchronous decentralized variants of \eqref{eq:sp_primal}-\eqref{eq:sp_dual}

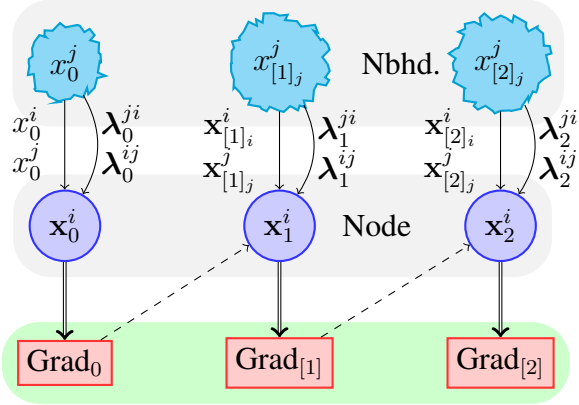
\begin{figure}
	\centering
	\begin{tikzpicture}[font=\large] 
	
	\matrix[row sep=1cm,column sep=1.5cm, ampersand replacement=\&] {
		\node (n_t-1) [neighbor] {$x^{j}_0$}; \&
		\node (n_t)   [neighbor] {$x^{j}_{[1]_j}$};     \&
		\node (n_t+1) [neighbor] {$x^{j}_{[2]_j}$};         \&

		\\
		\node (x_t-1) [agent] {$\bbx^{i}_0$}; \&
		\node (x_t)   [agent] {$\bbx^{i}_1$};     \&
		\node (x_t+1) [agent] {$\bbx^{i}_2$};         \&  
		\\
		
		\node (f_t-1) [local_nat]  {\text{Grad}$_{0}$};       \&
		\node (f_t)   [local_nat] {\text{Grad}$_{[1]}$};        \&
		\node (f_t+1) [local_nat] {\text{Grad}$_{[2]}$};       \&
		\\
		%
	};
	
	\path[->]
	(n_t-1) edge  [bend left] node [midway,right]  {\hspace{-2.55mm}$\begin{array}{c}
		\bblambda^{ji}_0\\
		\bblambda^{ij}_0
		\end{array}$} (x_t-1) 
	(n_t) edge  [bend left] node[midway, right]  {\hspace{-2.55mm}$\begin{array}{c}
		\bblambda^{ji}_1\\
		\bblambda^{ij}_1
		\end{array}$} (x_t)
	(n_t+1) edge  [bend left] node  [midway, right]{\hspace{-2.55mm}$\begin{array}{c}
		\bblambda^{ji}_2\\
		\bblambda^{ij}_2
		\end{array}$} (x_t+1)
	
	
	(n_t-1) edge  node [midway, left]  {$\begin{array}{c}
		x^{i}_0 \\
		x^{j}_0 
		\end{array}$} (x_t-1)
	(n_t) edge node [midway,left] {$\begin{array}{c}
		\bbx^{i}_{[1]_i} \\
		\bbx^{j}_{[1]_j}
		\end{array}$} (x_t)
	(n_t+1) edge node [midway, left] {$\begin{array}{c}
		\bbx^{i}_{[2]_i} \\
		\bbx^{j}_{[2]_j}
		\end{array}$} (x_t+1)

	(x_t-1) edge[double](f_t-1)				
	(x_t)   edge[double] (f_t)
	(x_t+1) edge [double](f_t+1)

	(f_t-1) edge [ dashed]  (x_t)				
	(f_t)   edge[ dashed]   (x_t+1)

	
	;
	
	\begin{pgfonlayer}{background}
	\node [background,
	fit=(n_t-1) (n_t+1), label=center:\large \hspace{30mm}Nbhd.] {};
	
	\node [background,
	fit=(x_t-1) (x_t+1), label=center:\large \hspace{25mm}Node] {};
	
	\node [background2,
	fit=(f_t-1) (f_t+1)] {};
	
	\end{pgfonlayer}
	\end{tikzpicture}
	\caption{\colb{Message passing in proposed algorithm. The red block represents the gradients required as mentioned in step 4 of the Algorithm \ref{algo_1} which depends upon the realizations $\bbtheta_{[t]_i}^i$ and $\bbtheta_{[t]_j}^j$ at time step $t$.}}\label{cartoon}
	\vspace{-8mm}
\end{figure}

In particular, to ameliorate the computational bottleneck associated with synchronized computation and communication rounds among the nodes, we consider situations in which observations and updates are subject to stochastic delays, i.e., an {\emph{asynchronous processing architecture}}. These delays take the form of random delays on the gradients which are used for the algorithm updates. We associate to each node $i$ in the network a time-dependent delay $\tau_i(t)$ for its stochastic gradient. Since the gradient corresponding to node $i$ are delayed by $\tau_i(t)$, it implies that the received gradient corresponds to $t-\tau_i(t)$ time slot which we denote as ${[t]_i}$. Rather than waiting for the current gradient at time $t$, agent $i$ instead uses the delayed gradient from the neighboring nodes at time $[t]_j$ for its update at time $[t]_i$. This leads to the following asynchronous primal update for stochastic online saddle point algorithm at each node $i$  
\begin{align}\label{eq:local_primal_update_asyn}
\bbx^{i}_{t+1}  =& \ccalP_{\ccalX}\Big[  \bbx^{i}_t  - \eps  
\Big(   \nabla_{\bbx^i} f^i (\bbx^{i}_{[t]_i} , \bbtheta^{i}_{[t]_i} )  \\
&+\!\!\sum_{j \in n_i}\!\! \left(\lam^{ij}_{t}+\lam^{ji}_{t}\right) \nabla_{\bbx^i} h^{ij}  \left(\bbx^{i}_{[t]_i}, \bbx^{j}_{[t]_j},\bbtheta^{i}_{[t]_i},\bbtheta^{j}_{[t]_j}\right)  \Big)\Big].\nonumber
\end{align}
Likewise, the dual update for each edge $(i,j)\in\mathcal{E}$ is
\begin{equation} \label{eq:local_dual_update_asyn}
\lam^{ij}_{t+1} \!\!=  \! \Big[\! (1\! - \eps^2 \delta) \lam^{ij}_{t}\! + \eps \left(\!h^{ij}\!\!\left(\!\bbx^{i}_{[t]_i}, \bbx^{j}_{[t]_j},\bbtheta^{i}_{[t]_i}\!,\!\bbtheta^{j}_{[t]_j}\!\right)\!\right)\!\!\Big]_+  .
\end{equation}
Note that to perform the asynchronous primal updates at node $i$ in \eqref{eq:local_primal_update_asyn}, delayed primal gradients $\nabla_{\bbx^i} f^i (\bbx^{i}_{[t]_i} , \bbtheta^{i}_{[t]_i} )$ and $\nabla_{\bbx^i} h^{ij}  \left(\bbx^{i}_{[t]_i}, \bbx^{j}_{[t]_j},\bbtheta^{i}_{[t]_i},\bbtheta^{j}_{[t]_j}\right)$ are utilized. Similarity, the dual delayed gradient is utilized for the update in \eqref{eq:local_dual_update_asyn}.  For the consistency in the algorithm implementation, it is assumed that at each node $i$, only the recent received copy of the gradient is kept and used for the update. Equivalently, this condition can be mentioned as 
%
${[t]_i}\geq [t-1]_i $
%
which implies that $\tau_i(t)\leq \tau_i(t-1)+1$. For brevity, we will use the notation $[\t]$ as a collective notation for all the delayed time instances as
%
$[\t]:=[[t]_1 ; \cdots ; [t]_N]$. 
%
The asynchronous algorithm is summarized in Algorithm~\ref{algo_1}. \colb{The practical implementation of the proposed asynchronous algorithm is explanied with the help of diagram in Fig.\ref{cartoon}. As described in figure, each node receives delayed parameters, gradients and carries out the updated accordingly.} The convergence guarantees for the proposed algorithm are shown to hold as $t-[t]_i\leq \tau$ is finite (assumption \ref{A6}) for all $i$ and $t$. The convergence results presented in \cite{koppel2017proximity} can be obtained as a special case with $\tau_{i}(t)=0$ from the results developed in this paper. This shows the generalization of the existing results in literature.  

Before proceeding with the convergence analysis, in order to have a tractable derivation, let us define the compact notation for the primal and dual delayed gradient as follows 
\begin{align}\label{d11}
\nabla_{\x^i}\hat{\mathcal{L}}^i_{[\t]}&\Big(\bbx^{i}_{[t]_i},\bbx^{j}_{[t]_j},\bblam^i_{[t]_i},\bblam^j_{[t]_j} \Big)
{:=}   \Big(   \nabla_{\bbx^i} f^i (\bbx^{i}_{[t]_i} , \bbtheta^{i}_{[t]_i} )  \\
&\ \ +\!\!\sum_{j \in n_i}\!\! (\lam^{ij}_{t}+\lam^{ji}_{t}) \nabla_{\bbx^i} h^{ij}  \left(\bbx^{i}_{[t]_i}, \bbx^{j}_{[t]_j},\bbtheta^{i}_{[t]_i},\!\bbtheta^{j}_{[t]_j}\right)  \! \!\Big) \nonumber
\end{align}
which follows from \eqref{eq:local_primal_update_asyn} and
\begin{align}\label{d21}
&\!\!\nabla_{\bblam^i}\hat{\mathcal{L}}^i_{[\t]}\left(\bbx^{i}_{[t]_i},\bbx^{j}_{[t]_j},\bblam^i_{[t]_i} \! \right){:=}h^{ij}\left(\!\bbx^{i}_{[t]_i}, \bbx^{j}_{[t]_j},\bbtheta^{i}_{[t]_i},\bbtheta^{j}_{[t]_j}\!\right)
\end{align} 
follows from \eqref{eq:local_dual_update_asyn}. These definitions of \eqref{d11} and \eqref{d21} will be used in rest of the places in this paper.

\begin{algorithm}[t] 
	\caption{ASSP: Asynchronous Stochastic Saddle Point}
	{\small \label{algo_1} 
		\begin{algorithmic}[1]
			\Require initialization $\bbx_{0}$ and $\bblam_{0}=\bb0$, step-size $\eps$, regularizer $\delta$
			\For {$t=1,2,\ldots, T$}
			\Loop {{\bf{ in parallel}} agent $i\in \mathcal{V}$}
			\State Send dual vars. $\bblam_{ij,t}$ to nbhd. $j\in n_i$
			\State Observe delayed gradients $\nabla_{\bbx^i} f^i (\bbx^{i}_{[t]_i} , \bbtheta^{i}_{[t]_i} )$, $\nabla_{\bbx^i} h^{ij}  \left(\bbx^{i}_{[t]_i}, \bbx^{j}_{[t]_j},\bbtheta^{i}_{[t]_i},\bbtheta^{j}_{[t]_j}\right)$ and constraint function $h^{ij}\left(\bbx^{i}_{[t]_i}, \bbx^{j}_{[t]_j},\bbtheta^{i}_{[t]_i},\bbtheta^{j}_{[t]_j}\right)$.
			\State Update $\bbx^{i}_{t+1}$ using \eqref{eq:local_primal_update_asyn}
			local parameter $\bbx^{i}_{t}$			\vspace{-1mm}
			\begin{align*}
			\bbx^{i}_{t+1} 
			=& \ccalP_{\ccalX}\Big[  \bbx^{i}_t  - \eps  
			\Big(   \nabla_{\bbx^i} f^i (\bbx^{i}_{[t]_i} , \bbtheta^{i}_{[t]_i} ) \nonumber
			\\
			&+\!\!\sum_{j \in n_i}\!\! \left(\lam^{ij}_{t}+\lam^{ji}_{t}\right) \nabla_{\bbx^i} h^{ij}  \left(\bbx^{i}_{[t]_i}, \bbx^{j}_{[t]_j},\bbtheta^{i}_{[t]_i},\bbtheta^{j}_{[t]_j}\right)  \Big)\Big]
			\end{align*}
			\vspace{-2mm}
			%
			%
			\State Update dual variables at each agent $i$ [cf. \eqref{eq:local_dual_update_asyn}]			\vspace{-1mm}
			\begin{align*} 
			\lam^{ij}_{t+1} \!=  \! \Big[\! (1\! - \eps^2 \delta) \lam^{ij}_{t} + \eps \left(h^{ij}\left(\bbx^{i}_{[t]_i}, \bbx^{j}_{[t]_j},\bbtheta^{i}_{[t]_i},\bbtheta^{j}_{[t]_j}\right)\!\right)\!\!\Big]_+  .
			\end{align*}			\vspace{-3mm}
			\EndLoop
			\EndFor
		\end{algorithmic}\normalfont}
\end{algorithm}

\section{Convergence in Expectation}\label{sec:convergence}
In this section, we establish convergence in expectation of the proposed asynchronous technique in \eqref{eq:local_primal_update_asyn}-\eqref{eq:local_dual_update_asyn} to a primal-dual optimal pair of the  problem formulated in \eqref{eq:coop_stoch_opt} when constant step sizes are used. Specifically, a sublinear bound on the average objective function optimality gap $F(\bbx_t) - F(\bbx^*)$ and the network-aggregate delayed constraint violation is established, both on average. The optimal feasible vector $\bbx^*$ is defined by \eqref{eq:coop_stoch_opt}. It is shown that the time-average primal objective function $F(\bbx_t)$ converges to the optimal value $F(\x^\star)$ at a rate of $\ccalO(1/\sqrt{T})$. Similarly, the time-average aggregated delayed constraint violation over network vanishes with the order of $\ccalO(T^{-1/4})$, both in expectation, where $T$ is the final iteration index.
To prove convergence of the stochastic asynchronous  saddle point method, some assumptions related to the system model and parameters are required which we state as follows.
\vspace{-3mm}
%
\begin{assumption}\label{A1}
	(\emph{Network connectivity}) The network $\ccalG$ is symmetric and connected with diameter $D$.
\end{assumption}
\vspace{-5mm}

\begin{assumption}\label{A2}
	(\emph{Existence of Optima})
	The set of primal-dual optimal pairs $\ccalX^* \times \bbLam^*$  of the constrained problem \eqref{eq:coop_stoch_opt} has non-empty intersection with the feasible domain $\ccalX^N \times \reals^M_+$.
\end{assumption}
\vspace{-5mm}
\begin{assumption}\label{A3}
	(\emph{Stochastic Gradient Variance}) The instantaneous objective and constraints for all $i$ and $t$ satisfy
	\begin{align}\label{as:third} 
	&\mathbb{E}\norm{ \nabla_{\bbx^i} f^i(\bbx^i,\bbtheta^{i}_{t})}^2\leq \sigma_{f}^2
	\\
	&\mathbb{E}\norm{\nabla_{\bbx^i} h^{ij}\left(\bbx^{i}, \bbx^{j},\bbtheta^{i}_{[t]_i},\bbtheta^{j}_{[t]_j}\right) }^2\leq\sigma_{h}^2
	\end{align}
	which states that the second moment of the norm of objective and constraint function gradients are bounded above.
\end{assumption}
\vspace{-5mm}
\begin{assumption}\label{A4} (\emph{Constraint Function Variance}) For the instantaneous constrain function for all pairs $(i,j)\in\mathcal{E}$ and $t$ over the compact set $\mathcal{X}$, it holds that 
	\begin{align} \label{as:fourth}
	\max_{(\bbx^{i},\bbx^{j})
		\in\mathcal{X}} \Ex{\left(h^{ij}\left(\bbx^{i}, \bbx^{j},\bbtheta^{i}_t,\bbtheta^{j}_t\right)^2\right)}\leq \sigma_{\bblam}^2
	\end{align}
	which implies that the maximum value the constraint function can take is bounded by some finite scalar $\sigma_{\bblam}^2$ in expectation.
\end{assumption}
\vspace{-5mm}
\begin{assumption}\label{A5}
	(Lipschitz continuity) The expected objective function defined in \eqref{eq:dist_stoch_opt} satisfies 
	\begin{align}
	\norm{F(\x)-F(\y)}\leq L_f\norm{\x-\y}.
	\end{align}
	for any $(\x,\y)\in \mathbb{R}^{Np}$.
\end{assumption}
\vspace{-5mm}
\begin{assumption}\label{A6}
	(Bounded Delay) The delay  $\tau_i(t)$ associated with each node $i$ is upper bounded: $\tau_i(t)\leq\tau$ for some $\tau<\infty$.
\end{assumption}
Assumption~\ref{A1} ensures that the graph is connected and the rate at which information diffuses across the network is finite. This condition is standard in distributed algorithms \cite{RamNedicVeeravalli,1506308}. 
Moreover, Assumption~\ref{A2} \colb{is a Slater's condition which} makes sure the existence of an optimal primal-dual pair within the feasible sets onto which projections occur which are necessary for various quantities to be bounded. It has appeared in various forms to guarantee existence of solutions in constrained settings \cite{bertsekas03}. 
Assumption \ref{A3} assumes an upper bound on the mean norm of the primal and dual stochastic gradients, which is crucial to developing the gradient bounds for the Lagrangian used in the proof. Assumption \ref{A4} yields an upper  bound on the maximum possible value of the constraint function in expectation similar to that of \cite{mahdavi2012trading}, and is guaranteed to hold when $\ccalX$ is compact and $h^{ij}$ is Lipschitz. Assumption ~\ref{A5} is related to the Lipschitz continuity of the primal objective function. Assumption ~\ref{A6} ensures that the delay is always bounded by $\tau$, which holds in most wireless communications problems and autonomous multi-agent networks \cite{nedic2001distributed}. 

For the analysis to follow, we first derive bounds on the mean square-norms of the stochastic gradients of the Lagrangian. Thus, consider the mean square-norm of the primal stochastic gradient of the Lagrangian, stated as:
\begin{align}\label{main}
\!\!\mathbb{E}[\! \|\nabla_{\x}\hat{\mathcal{L}}_{t}(\x,\bblam) \|^2]\!
&\leq \! {2}N\!\left[\max_i\mathbb{E}\norm{ \nabla_{\bbx^i} f^i(\bbx^i,\bbtheta^{i}_{t})}^2\right]
\\
& +{2M^2}\!\!\norm{\boldsymbol{\lambda}}^2\!\!\max_{(i,j)\in\mathcal{E}} \mathbb{E}\!\norm{\nabla_{\bbx^i} h^{ij} \!\!\left(\bbx^{i}, \bbx^{j},\bbtheta^{i}_{t},\bbtheta^{j}_{t}\!\right)\!}^2\nonumber
\end{align}
Now apply Assumptions in \ref{A3} to the mean square-norm terms in \eqref{main} to obtain 
\begin{align}
\mathbb{E}[ \|\nabla_{\x}\hat{\mathcal{L}}_{t}(\x,\bblam) \|^2]
&\leq {2}N\sigma_{f}^2+{2M^2}\norm{\boldsymbol{\lambda}}^2\sigma_{h}^2\nonumber \\
&\leq {2}{(N+M^2)} L^2 (1+\norm{\boldsymbol{\lambda}}^2)\label{gradient_1}
\end{align}
where, $L^2=\max(\sigma_{f}^2,\sigma_{h}^2)$. Similarly for the gradient with respect to the dual variable $\bblam$, we have  
\begin{align}
\!\E{ \! \| \nabla_{\bblam}\hat{\mathcal{L}}_{t}(\x,\bblam\!) \|^2\! }\!
&\!\leq
\! 2M \!\!\!\max_{(i,j)\in\mathcal{E}} \! \! \!\mathbb{E}[\!(h^{ij}\!(\bbx^{i}, \!\bbx^{j},\!\bbtheta^{i}_t,\!\bbtheta^{j}_t) \!)^2\!]
\! \! +\!  2\delta^2\! \eps^2 \! \| \bblambda \|^2   \label{gradient_2}\nonumber \\
&\leq 
\! 2M\sigma_{\boldsymbol{\lambda}}^2+\!  2\delta^2\! \eps^2 \| \bblambda \|^2   
\end{align}
The bounds developed in \eqref{gradient_1} and \eqref{gradient_2} are in terms of the norm of the dual variable and utilizes the Assumptions \ref{A3} and \ref{A4}. It is important to note that these bounds are for arbitrary $\x$ and $\bblam$, therefore holds for any realization of the primal $\x_t$ and dual variables $\bblam_t$. Before proceeding towards the main lemmas and theorem of this work, a remark on the importance of the bounds in \eqref{gradient_1} and \eqref{gradient_2} is due.  

\begin{remark} \normalfont
	Conventionally, in the analysis of primal-dual methods, the primal and dual gradients are bounded by constants. In contrast, here our upper-estimates depend upon the magnitude of the dual variable $\bblam$ which allows us to avoid dual set projections onto a compact set, and instead operate with unbounded dual sets $\reals_{+}^M$. We are able to do so via exploitation of the dual regularization term $-(\delta \eps/2)\|\bblambda\|^2$ that allows us to control the growth of the constraint violation.
\end{remark}
%
Subsequently, we establish a lemma for the instantaneous Lagrangian difference $\hat{\mathcal{L}}_{[\t]}(\x_{[\t]},\bblam) -\hat{\mathcal{L}}_{[\t]}(\x, \bblam_{t}) $ by a telescopic quantity involving the primal and dual iterates, as well as the magnitude of the primal and dual gradients. This lemma is crucial to the proof of our main result at the end of this section.
%
\begin{lemma}\label{lemma1}
	Under the Assumptions ~\ref{A1} - ~\ref{A6} ,  the sequence $(\bbx_t, \bblam_t )$ generated by the proposed asynchronous stochastic saddle point algorithm in \eqref{eq:local_primal_update_asyn}-\eqref{eq:local_dual_update_asyn} is such that for a constant step size $\eps$, the  instantaneous Lagrangian difference sequence $\hat{\mathcal{L}}_{[\t]}(\x_{[\t]},\bblam) -\hat{\mathcal{L}}_{[\t]}(\x, \bblam_{t})$ satisfies the decrement property
	\begin{align}\label{eq:lemma12}
	\hat{\mathcal{L}}_{[\t]}&(\x_{[\t]},\bblam) -\hat{\mathcal{L}}_{[\t]}(\x, \bblam_{t})  \nonumber 
	\\ 
	\leq &\frac{1}{2\eps} \big( \|\bbx_{t}\!-\!\bbx\|^2 \!- \!\|\bbx_{t+1}\!-\!\bbx\|^2 
	+\|\bblam_{t}-\bblam\|^2 \!-\! \|\bblam_{t+1}\!-\!\bblam\|^2 \big)  \nonumber
	\\
	&+\frac{\eps}{2} \left(\norm{\nabla_{\bblam}\hat{\mathcal{L}}_{[\t]}(\x_{[\t]},\bblam_{t})}^2 +\| \nabla_{\bbx}\hat{\mathcal{L}}_{[\t]}(\x_{[\t]},\bblam_{t})\|^2\right)\nonumber
	\\
	&\quad+\ip{\nabla_{\bbx}\hat{\mathcal{L}}_{[\t]}(\x_{[\t]},\bblam_{t}), (\bbx_{[\t]}-\bbx_{t})}
	\end{align}
\end{lemma}

\begin{myproof}
	See Appendix A.
\end{myproof}

Lemma \ref{lemma1} exploits the fact that the stochastic augmented Lagrangian is convex-concave with respect to its primal and dual variables to obtain an upper bound for the difference $\hat{\mathcal{L}}_{[\t]}(\x_{[\t]},\bblam) -\hat{\mathcal{L}}_{[\t]}(\x, \bblam_{t})$ in terms of the difference between the current and the next  primal and dual iterates to a fixed primal-dual pair $(\bbx, \bblam)$, as well as the square magnitudes of the primal and dual gradients. Observe that here, relative to \cite{koppel2017proximity}[Proposition 1], an additional term is present which appears that represents the directional error caused by asynchronous updates of Algorithm \ref{algo_1}.
This contractive property is the basis for establishing the convergence of the primal iterates to their constrained optimum given by \eqref{eq:coop_stoch_opt} in terms of mean objective function evaluation and mean constraint violation with constant step-size selection. Before proceeding towards  our main theorem, we establish an additional lemma which simplifies its proof and clarifies ideas.
\begin{lemma}\label{lemma12}
	Denote as $(\bbx_t, \bblam_t )$ the sequence generated by the asynchronous saddle point algorithm in \eqref{eq:local_primal_update_asyn}-\eqref{eq:local_dual_update_asyn} with stepsize $\eps$. If Assumptions \ref{A1} - \ref{A6} holds,  then it holds that 
	\begin{align}\label{eq:lemma122}
	\mathbb{E}\Big[&\sum_{t=1}^T [F(\bbx_t) \!- \! F(\bbx)] 
	\!+\!\!\!\!\sum_{(i,j)\in \mathcal{E}}\!\!\!\!
	\lam^{ij}\!\!\left(\sum_{t=1}^Th^{ij}\!\!\left(\!\bbx^{i}_{[t]_i},\!\bbx^{j}_{[t]_j},\!\bbtheta^i_{[t]_i},\!\bbtheta^j_{[t]_j} \!\right)\!\!\right) \nonumber \\
	&- \Big(\!\frac{\delta \eps T}{2} \!+\! \frac{1}{2\eps}\Big)\! (\lambda^{ij})^2 \Big] 
	\leq    \frac{1}{2\eps} \|\bbx_1 \!-\! \!\bbx \|^2+  \frac{ \eps T K}{2}
	\end{align}
	where the constant $K$ is defined in terms of system parameters as 
	\begin{align}\label{cosntant}
	K:= M\sigma^2_{\bblam} +(N+M)L[(1/2)L+\tau^2(2L+L_f)].
	\end{align}
\end{lemma}

\begin{myproof}
	See Appendix B.
\end{myproof}


Lemma \ref{lemma12} is derived by considering Lemma~\ref{lemma1}, computing expectations, and applying  \eqref{gradient_1} and \eqref{gradient_2}. This lemma describes the global behavior of the augmented Lagrangian when following Algorithm \ref{algo_1}, and may be used to establish convergence in terms of primal objective optimality gap and aggregated network constraint violation, as we state in the following theorem.
%
\begin{thm}\label{theorem1}
	Under the Assumptions \ref{A1} - \ref{A6} , denote $(\bbx_t, \bblam_t )$ as the sequence of primal-dual variables generated by Algorithm \ref{algo_1} [cf. \eqref{eq:local_primal_update_asyn}-\eqref{eq:local_dual_update_asyn}]. When the algorithm is run for $T$ total iterations with constant step size $\eps=1/\sqrt{T}$,  the average time aggregation of the sub-optimality sequence $\Ex{{F(\bbx_t) - F(\bbx^*)}}$, with $\bbx^*$ defined by \eqref{eq:coop_stoch_opt}, grows sublinearly with $T$ as 
	\begin{align}\label{thm1}
	\sum_{t=1}^T \mathbb{E}[F(\bbx_t) \!- \! F(\bbx^*)] 
	\!& \leq    \ccalO(\sqrt{T}).
	\end{align}
	Likewise, the delayed time aggregation of the average constrain violation also grows sublinearly in $T$ as
	\begin{align} \label{eq:theorem12}
	&\hspace{-6mm}\!\!\sum_{(i,j)\in\mathcal{E}}\!\!\!\mathbb{E}\Bigg[\!\!\Big[\!\!\sum_{t=1}^T\! \left(h^{ij}(\bbx^{i}_{[t]_i},\bbx^{j}_{[t]_j},\bbtheta^i_{[t]_i},\bbtheta^j_{[t]_j})-\gamma_{ij}\right) \Big]_+\!  \Bigg]\leq \ccalO(T^{3/4}).
	%
	\end{align}
\end{thm}
%
\begin{myproof} 
	See Appendix C.
\end{myproof} 

Theorem \ref{theorem1} presents the behavior of the Algorithm \ref{algo_1} when run for $T$ total iterations with a constant step size. Specifically, the average aggregated objective function error sequence is upper bounded by a constant time $\sqrt{T}$ sequence. This establishes that the expected value of the objective function $\Ex{F(\x_t)}$ will become closer to the optimal $F(\x^\star)$ for larger $T$. Similar behavior is shown by the average delayed aggregated network constraint violation term. The key innovation establishing this result is the bound on the directional error caused by asynchrony. We achieve this by bounding it in terms of the gradient norms and dual variable $\bblam$ and exploiting the fact that the delay is at-worst bounded.

These results are similar to those for the unconstrained convex optimization problems with sub-gradient descent approach and constant step size. For most of the algorithms in this context \cite[Section 2.2, eqn. 2.19]{nemirovski2009robust}, or \cite[Section 4]{bach2014adaptivity}, convergence to the neighborhood of size $\mathcal{O}(\eps T )$ is well known. For such algorithms, the primal sub-optimality is shown of the order $\ccalO (\ep T)$ is shown and the radius of suboptimality is optimally controlled by selecting $\ep=1/\sqrt{T}$ \cite{Zinkevich2003}. The bound $\ccalO(T^{3/4})$ on the constraint violation aggregation is comparable to existing results for synchronized multi-agent online learning \cite{mahdavi2012trading} and stochastic approximation \cite{koppel2017proximity}.

Moreover, tt is possible to extend the result in \eqref{thm1} to show that the convergence results of the average objective function error sequence also holds for the running average of primal iterates $\hat{x}_T:=\frac{1}{T}\sum_{t=1}^{T}\x_t$ as in \cite[corollary 1]{koppel2017proximity}. However,  obtaining such a result for the constraint violation in \eqref{eq:theorem12} is not straightforward due to the presence of delayed primal variables.  Subsequently, we turn to studying the empirical performance of Algorithm \ref{algo_1} for developing intelligent interference management in communication systems.

\section{Interference management through pricing}\label{sec:fields}
The rising number of cellular users has fueled the increase in infrastructure spending by cellular operators towards better serving densely populated areas. In order to circumvent the near-absolute limits on spectrum availability, the current and future generations rely heavily on frequency reuse via small cells and associated interference management techniques \cite{haykin2005cognitive,lopez2009ofdma,base}. This work builds upon the pricing-based interference management framework proposed in \cite{base}. We consider heterogeneous networks with multiple autonomous small cell users. Under heavy load situations, the macro base station (MBS) may assign the same operating frequency to multiple but geographically disparate small cell base stations (SCBS) and macro cell users (MU). The base station regulates the resulting cross-tier interference (from SCBS to MU) by penalizing the received interference power at the MUs. Consequently, the SCBSs coordinate among themselves and employ power control to limit their interference at the MUs. This section considers the pricing problem from the perspective of the BS that seeks to maximize its revenue. 
\begin{figure}
	\vspace{-8mm}
	\centering
	\scriptsize
	\includegraphics[scale=0.25]{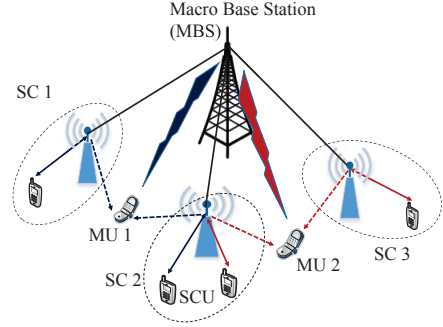}
	\caption{Heterogeneous cellular network with one MBS, two MUs, and three SCBSs with each serving one, two, and one SCU, respectively.  }
	\label{model}
	\vspace{-5mm}
\end{figure}

\subsection{Problem formulation}
Consider the network depicted in Fig. \ref{model}, consisting of a MBS serving $M$ MU users and $N$ SCBSs \cite{base}. Each MU is assigned a unique subchannel, indexed by $i \in \{1, \ldots, M\}$. At times of high traffic, the BS also allows the SCBSs to use these $M$ subchannels, so that the $n$-th SCBS may serve $K_n \leq M$ SCUs. In other words, at each time slot, a particular subchannel is used by MU $i$ and a non-empty set of SCBSs $\mathcal{N}_i \subset \{1, \ldots, N\}$. Denoting the channel gains between the $n$-th SCBS and $i$-th MU by $g_{ni}$, it follows that the total interference at the MU $i$ is given by
$I_i := \sum_{n \in \mathcal{N}_i} g_{ni} p_n^i$
where $p_n^i$ is the transmit power of SCBS $n$ while using subchannel assigned to MU $i$. The BS regulates this cross-tier interference by imposing a penalty $x_n^i$ on the SCBSs $n \in \mathcal{N}_i$. The total revenue generated by the BS is therefore given by
\begin{align}
\sum_{i=1}^M \sum_{n \in \mathcal{N}_i} x_n^i g_{ni} p_n^i
\end{align}
which the BS seeks to maximize. The BS also adheres to the constraint that the total penalty imposed on each SCBS is within certain limit, i.e.,
\begin{align}
C_{\min} \leq \sum_{i: n \in \mathcal{N}_i} x_n^i \leq C_{\max}.
\end{align}
The limit on the maximum and minimum penalties can also be viewed as a means for BS to be fair to all small cell operators.
The power allocation at the SCBSs is governed by their local transmission costs, denoted by $c$ per unit transmit power, and the interference prices levied by the BS. As in \cite{base}, each SCBS solves a penalized rate minimization subproblem, resulting in the power allocation
\begin{align}\label{power_func}
\vspace{-0mm}
p_n^i = \left((W/(c\mu_n+ \nu_n x_n^i)) - (1/h_n^i)\right)_{+}
\end{align}
for all $n \in \mathcal{N}_i$ and $1\leq i\leq M$. Here, $h_n^i$ is the channel gain between $n$-th SCBS and its scheduled user,  $\mu_n$ and $\nu_n$ represent SCBS-specific parameters used to trade-off the achieved sum rate against the transmission costs and $W$ is the bandwidth per subcarrier. Finally, the channel gains $g_{ni}$ and $h_n^i$ are not known in advance, so the BS seeks to solve the following stochastic optimization problem:
\begin{subequations}\label{prob_for_iain}
	\begin{align}
	\max_{\{x_n^i\}} &\sum_{i=1}^M \sum_{n \in \mathcal{N}_i} \Ex{x_n^i g_{ni} p_n^i (x_n^i, h_n^i)}\label{obj_num} 
	\\
	\text{s. t. } &\sum_{n \in \mathcal{N}_i} \Ex{g_{ni} p_n^i (x_n^i, h_n^i)} \leq \gamma_{i}  ~ 1\leq i \leq M \label{const1_num}
	\\
	&C_{\min} \leq \sum_{i:n\in\mathcal{N}_i} x_n^i \leq C_{\max} ~~1\leq n \leq N\label{const2_num}
	\end{align}
\end{subequations}
Here, $\gamma_i$ is the interference power margin and observe that the interference constraint is required to hold only on an average, while the limits on the interference penalties are imposed at every time slot. \blue{It is remarked that the similar pricing based interference management scheme is considered is \cite{base} but with the assumption that the distribution of the random variables is known. For this work, we omit such assumption and propose an online solution to stochastic optimization problem in \eqref{prob_for_iain}.}
\begin{algorithm}
	\caption{Online interference management through pricing}
	{\small \label{algo_2} 
		\begin{algorithmic}[1]
			\Require initialization $\bbx_{0}$ and $\bblam_{0}=\bb0$, step-size $\eps$, regularizer $\delta$
			\For {$t=1,2,\ldots, T$}
			\Loop {{\bf{ in parallel}} for all MU and SCBS user}
			\State Send dual vars. $\lambda^{m}_{t}$ to nbhd.
			\State Observe the delayed primal and dual (sub)-gradients 
			\State Update the price $x^i_{n,{t+1}}$ at SCBS $n$ as 
			\begin{align*}
			\!\!\!\!\!\!\!{x}^i_{n,t+1}=&\ccalP_{\mathcal{X}_n}\!\!\Big[\!x^i_{n,t}\!+\ep\Big(g_{ni,[t]}\left[\! \frac{W(c\mu_n\!+\!\nu_n\lambda^i_t)}{{(c\mu_n\!+\!\nu_n x^i_{n,{[t]}})}^2}\!-\!\frac{1}{h^i_{n,[t]}}\!\right]\!\!\cdot \!\mathbf{1}(x^i_{n,[t]})\!\Big)\!\Big]
			\end{align*}
			\State Update dual variables at each MU $i$ [cf. \eqref{eq:local_dual_update_asyn}]			\vspace{-1mm}
			\begin{align*} 
			\hspace{-0.1cm}\!\!\!\!\!\lambda^i_{t+1}=&\left[(1\!\!+\!\!\delta \ep^2)\lambda^i_t\!-\!\eps\!\left(\!\!\gamma_i\!\!-\!\!\sum\limits_{n\in\mathcal{N}_i}\!\!g_{ni,[t]}\left(\!\! \frac{W}{c\mu_n+\nu_n x^i_{n,{[t]}}}-\frac{1}{h^i_{n,[t]}}\!\!\right)_{\!\!\!+}\!\right)\!\!\!\right]_{\!\!+}.
			\end{align*}			\vspace{-3mm}
			\EndLoop
			\EndFor
		\end{algorithmic}\normalfont}
\end{algorithm}\vspace{-5mm}
\subsection{Solution using stochastic saddle point algorithm}
It can be seen that the stochastic optimization problem formulated in \eqref{prob_for_iain}
is of the form required in \eqref{eq:coop_stoch_opt} with $\mathcal{X}$ capturing the constraint in \eqref{const2_num}. 
Since the random variables $h_n^i$ and $g_{ni,t}$  have bounded moments, the assumptions in Section \ref{sec:convergence} can be readily verified. Further, the saddle point method may be applied for solving \eqref{prob_for_iain}. To do so, we use the preceding definition of the power function $p_n^i$ defined as in \eqref{power_func}, and associating dual variable $\lambda^i$ with the $i$-th constraint in \eqref{prob_for_iain}, the stochastic augmented Lagrangian is given by 
\begin{align}
\hat{\mathcal{L}_t}(\x,\bblam)=&\sum_{i=1}^{M}\sum\limits_{n\in\mathcal{N}_i}{x_n^i g_{ni,t} p_n^i (x_n^i, h_{n,t}^i)}\label{stochs_lag_num}
\\
&+\sum_{i=1}^{M}{\lambda^i[\gamma_{i}\!\!-\!\!\!\!\sum\limits_{n\in\mathcal{N}_i}{g_{ni,t} p_n^i (x_n^i, h_{n,t}^i)}]\!}-\frac{\delta \ep}{2}\norm{\bblam}^2\nonumber
\end{align}
where $\x$ collects the variables $\{x_n^i\}_{i=1,n\in n_i}^M$ and $\bblam$ collects the dual variables $\{\lambda^i\}_{i=1}^M$. 

The asynchronous saddle point method for pricing-based interference management in wireless systems then takes the form of Algorithm \ref{algo_2} with the modified projection defined as
\begin{align}
P_{\mathcal{X}_n}(\u):= \min_{y} &\norm{\y-\u}\nonumber \\
\text{s. t.}\ \  &C_{\min} \leq \ip{\mathbf{1},\y} \leq C_{\max}.
\end{align}
In order to get the primal and dual updates of step $5$ and $6$, note that  the subgradient  of the Lagrangian in \eqref{stochs_lag_num} with respect to primal variables is given by 
\begin{align}
\partial_{x_n^i}&\mathcal{L}_t(x_{n}^i,\lambda^i):=g_{ni,t}\left[ \frac{W(c\mu_n+\nu_n\lambda^i)}{{(c\mu_n+\nu_n x^i_{n})}^2}-\frac{1}{h^i_{n,t}}\right]\cdot\mathbf{1}(x^i_{n})\nonumber
\end{align}
and the gradient of  the Lagrangian with respect to the dual variable $\bblam
^i$ is given by
\begin{align}
\nabla_{\lambda^i}\mathcal{L}_t(x_{n}^i,\lambda^i)=&\gamma_i\!-\!\!\!\!\sum\limits_{n\in\mathcal{N}_i}\!\!\!g_{ni,t}\!\left(\!\! \frac{W}{c\mu_n+\nu_n x^i_{n,t}}\!-\!\frac{1}{h_{n,t}^i}\!\right)_+\!\!\!-\delta \ep \lambda^i.\nonumber
\end{align}
Observe here that the implementation in Algorithm 2 allows the primal updates to be carried out in a decentralized manner. On the other hand, the base station carries out the dual updates. Consequently both, the SCBSs and the BSs may utilize old price iterates $x^i_{n,[t]}$.

\begin{figure*}
	\centering
	\setcounter{subfigure}{0}
	\begin{subfigure}{0.8\columnwidth}
		\includegraphics[width=\linewidth, height = 0.6\linewidth]{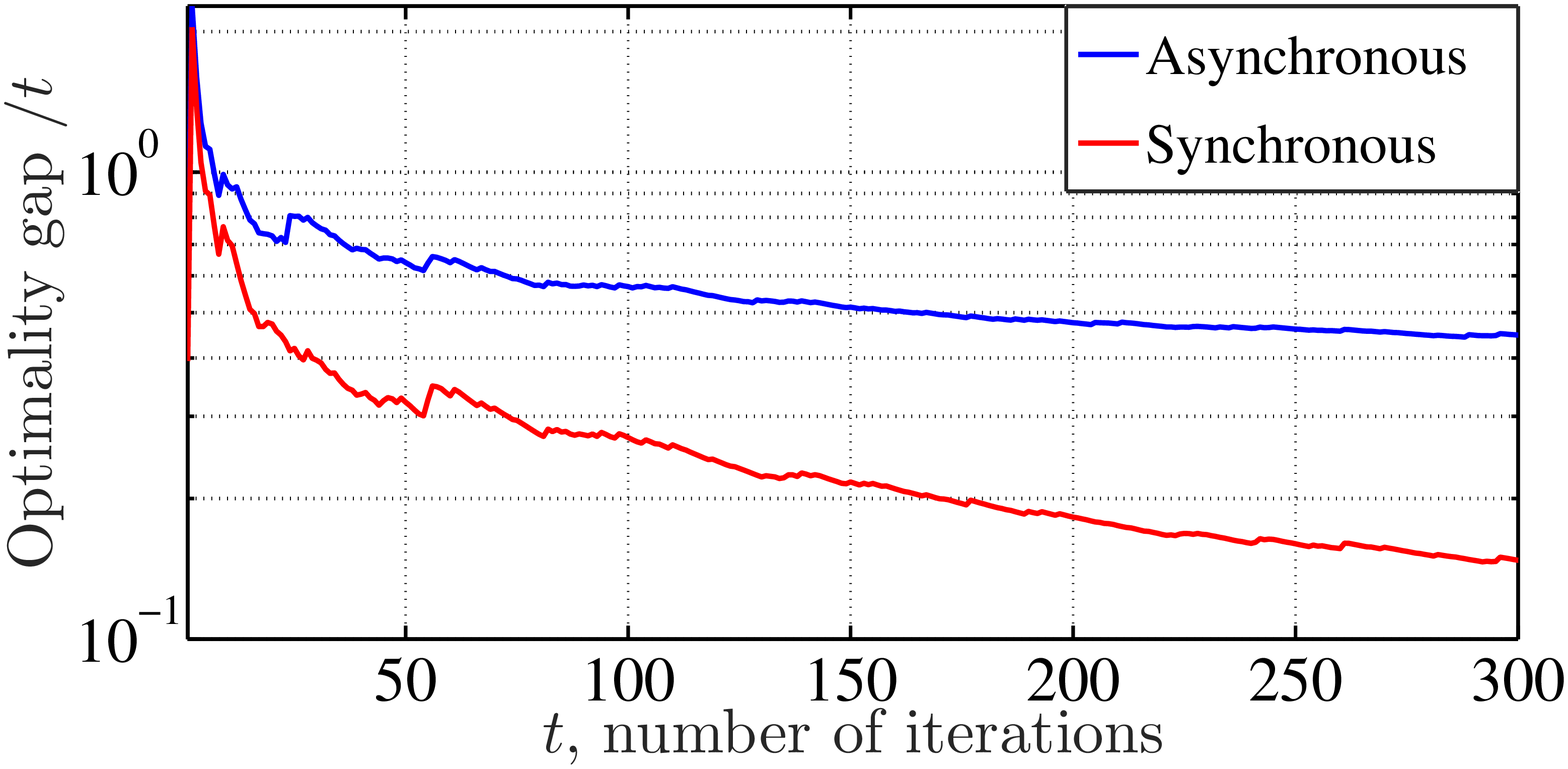}
		\caption{Average objective sub-optimality vs. iteration $t$}
		\label{obj}
	\end{subfigure}\hspace{2.8cm}
	\begin{subfigure}{0.8\columnwidth}
		\includegraphics[width=\linewidth, height = 0.6\linewidth]{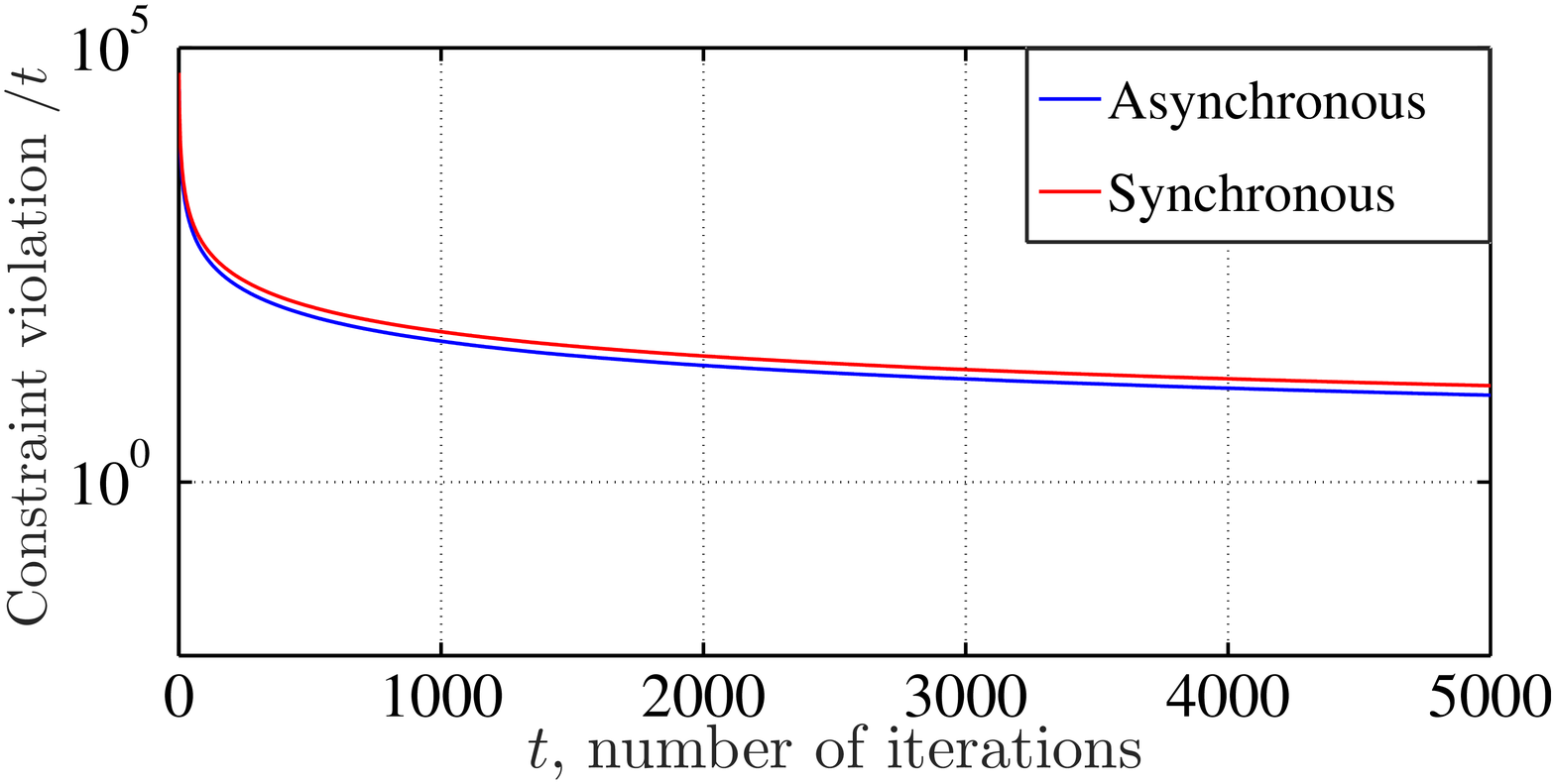}
		\caption{Average constraint violation vs. iteration $t$ }
		\label{vio}
	\end{subfigure}
	\caption{Algorithm \ref{algo_2} applied to a 5G cellular network with two MBS  user and three SCBSs. The $y$ axis of first figure is  $\frac{1}{t}\sum_{u=1}^t \mathbb{E}[F(\bbx_u) \!- \! F(\bbx^*)]$ and of second figure is $(1/t)\mathbb{E}\Big[\Big[\!\sum_{u=1}^t\! \Big(h^{i}\left(\{\bbx^{j},\bbtheta^{j}_{t}\}_{j\in n'_i}\right)\Big) \Big]_+\!  \Big]$.   Observe that both the asynchronous and synchronous implementations attain convergence but the asynchronous method settles to a higher level of sub-optimality. Thus, we may solve decentralized online learning problems without a synchronized clock.}
	\vspace{-6mm}
\end{figure*}
For the simulation purposes, we considered a cellular network with $M=2$ MBSs and $N=3$ SCBSs with index $\{m1,m2\}$ and $\{s1,s2,s3\}$. The scenario considered is similar to as shown in Fig.~\ref{model}, means that \{s1, s2\} are in the neighborhood of MU $m1$ and $\{s2,s3\}$ constitutes the neighborhood of MU $m2$.   The random channel gain $g_{ni}$ and $h_n^i$ are assumed to be exponentially distributed with mean $\mu=3$. The minimum and maximum values $C_{\min}=0.9$ and $C_{\max}=20$. The other parameter values are $W=1MHz$, $ \gamma_i=-3$ dB, $\delta=10^{-5}$, $c=0.1$, $\mu_n=\nu_n=1$, and $\epsilon=0.01$. The maximum delay parameter is $\tau=10$.

Fig.~\ref{obj} shows the  difference of running average of primal objective from its optimal value. It is important to note that the difference goes to zero as $t\rightarrow\infty$. The result for both synchronous and asynchronous algorithm algorithms are plotted. The optimal value to plot Fig.~\ref{obj} is obtained by running the synchronous algorithm for long duration of time and utilizing the converged value as the optimal one. We observe that running the saddle point method without synchrony breaks the bottleneck associated with heterogeneous computing capabilities of different nodes, although it attains slightly slower learning than its synchronized counterpart.
\begin{figure}
	\hspace{5mm}\includegraphics[width=0.8\linewidth, height = 0.4\linewidth]{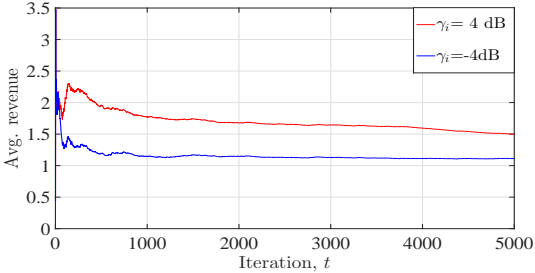}
	\caption{Average revenue for different interference power margin. }
	\label{revenue}
\end{figure}
Fig.~\ref{vio} shows the behavior of constraint violation term derived in \eqref{eq:theorem12} for a randomly chosen MBS. In the sample path of the empirical average constraint violation, the trend of sublinear growth of objective sub-optimality from Fig.~\ref{obj} is further substantiated by convergence in expectation of the constraint violation as the iteration index $t$ increases. We observe that the performance reduction of asynchronous operations relative to synchronous ones is smaller with respect to constraint violation as compared with primal-suboptimality, corroborating the rate analysis of Theorem \ref{theorem1}.

\begin{table}
	\begin{center}
		\begin{tabular}{ |c|c|c|c| } 
			\hline
			User & Algo. 2 & Naive approach \\
			\hline 
			MU $1$& $29$ dB & $22$ dB \\ 
			MU $2$& $28$ dB & $22$ dB \\ 
			\hline
		\end{tabular}
	\end{center}
	\caption{Comparison of Average SINR at MU}
	\label{table}
	\vspace{-7mm}
\end{table}

\colb{Fig.~\ref{revenue} shows that the average value of the revenue generated by  MBS converges to a higher value for higher interference power margin $\gamma_i=4$ dB. It shows the advantage of proposed interference management scheme which exploits the higher allowed interference as a resource to generate revenue for the MBS. Next, the proposed technique is compared with a 'naive' approach in which each SCBS user transmits at unity power all the times irrespective of the allowed interference power margin and channel conditions. The result in Table \ref{table} shows that the SINR achieved at both the MUs is higher for the proposed technique than that of the naive approach. It is due to the fact that proposed technique takes care of current channel conditions and therefore limits the interference caused to MU due to SCBS user transmission. }

%

\vspace{-3mm}
\section{Conclusion} \label{sec:conclusion}
We considered multi-agent stochastic optimization problems where the hypothesis that all agents are trying to learn common parameters may be violated, with the additional stipulation that agents do not even operate on a synchronized time-scale. To solve this problem such that agents give preference to locally observed information while incorporating the relevant information of others, we formulated this task as a decentralized stochastic program with convex proximity constraints which incentivize distinct nodes to make decisions which are close to one another. We derived an asynchronous stochastic variant of the Arrow-Hurwicz saddle point method to solve this problem through the use of a dual-augmented Lagrangian. 

We established that in expectation, under a constant step-size regime, the time-average suboptimality and constraint violation are contained in a neighborhood whose radius vanishes with increasing number of iterations (Theorem \ref{theorem1}). This result extends existing results for multi-agent convex stochastic programs with inequality constraints to asynchronous computing architectures which are important for wireless systems.

We then considered an empirical evaluation for the task of pricing-based interference management in large distributed wireless networks. For this application setting, we observe that the theoretical convergence rates translate into practice, and further that the use of asynchronous updates breaks the computational bottleneck associated with requiring devices to operate on a common clock. In particular, the asynchronous saddle point method learns slightly more slowly its synchronous counterpart while yielding a substantial complexity reduction on a real wireless application, and thus holds promise for other online multi-agent settings where synchronous operations and consensus are overly restrictive.

\vspace{-0mm}
\section*{Appendix 0: Derivation of Generalized Algorithm }\label{link}
Consider the generalized the problem in \eqref{eq:coop_stoch_opt_gen}. The stochastic augmented Lagrangian takes the form
\begin{align} \label{eq:stoch_lagrangian_gen}
\hat{\ccalL}_t(\bbx,\bblambda)&= \sum_{i=1}^N 
\bigg[f^i(\bbx^i,\bbtheta^{i}_t)\nonumber \\
&\ \ \ \ +\!\! \left(\!\!\!
\ip{\bblam^{i}\!, \h^{i}\left(\!\!\{\bbx^j,\bbtheta^{j}_t\}_{j\in n'_i} \right)\!\!}
-\frac{\delta \eps}{2} \norm{\bblam^{i}}^2\!\right)\!\!\bigg].
\end{align}
The primal update for the generalized asynchronous stochastic saddle point algorithm at each node $i$ is given by 
\begin{align}\label{eq:local_primal_update_asyn_gen}
\bbx^{i}_{t+1} 
= \ccalP_{\ccalX}\Big[  & \bbx^{i}_t  - \eps  
\Big(   \nabla_{\bbx^i} f^i (\bbx^{i}_{[t]_i} , \bbtheta^{i}_{[t]_i} ) \nonumber
\\
&\hspace{-3mm} + \!\!\sum_{k \in n'_i}\!\! \ip{\bblam^{k}_t, \nabla_{\bbx^i} \h_{k}  \left(\{\bbx^{j}_{[t]_j},\bbtheta^{j}_{[t]_j}\}_{j\in n'_k}\right)}  \Big)\Big].
\end{align}
through applying comparable logic to Section \ref{sec:algorithm}. Likewise, the dual variable updates at each node $i$ is given by	
\begin{equation} \label{eq:local_dual_update_asyn_gen}
\bblam^{i}_{t+1} \!=  \! \Big[\! (1\! - \eps^2 \delta) \bblam^{i}_t + \eps \left(\h^{i}\left(\{\bbx^{j}_{[t]_j},\bbtheta^j_{[t]_j}\}_{j\in n'_i} \right)\!\right)\Big]_+  .
\end{equation}
Before proceeding with the convergence analysis, to have a tractable derivation, let us define the compact notation for the primal and dual delayed gradient as follows 
\begin{align}\label{d1}
&\nabla_{\x^i}\hat{\mathcal{L}}^i_{[\t]}\left(\{\bbx^{j}_{[t]_j},\bblam^j_{[t]_j}\}_{j\in n'_i} \right)\boldsymbol{:=}   \Big(\nabla_{\bbx^i} f^i (\bbx^{i}_{[t]_i} , \bbtheta^{i}_{[t]_i} ) \nonumber
\\
&\hspace{16mm}+ \sum_{k \in n'_i} \ip{\bblam^{k}_t, \nabla_{\bbx^i} \h_{k}  \left(\{\bbx^{j}_{[t]_j},\bbtheta^{j}_{[t]_j}\}_{j\in n'_k}\right)}\Big)
\end{align} 
which generalizes \eqref{eq:local_primal_update_asyn} to the setting of \eqref{eq:coop_stoch_opt_gen} and the notation 
\begin{align}\label{d2}
&\!\!\nabla_{\bblam^i}\hat{\mathcal{L}}^i_{[\t]}\left(\{\bbx^{j}_{[t]_j},\bblam^j_{[t]_j}\}_{j\in n'_i} \right)\boldsymbol{:=}\h^{i}\left(\{\bbx^{j}_{[t]_j},\bbtheta^j_{[t]_j}\}_{j\in n'_i} \right)
\end{align} 
follows from \eqref{eq:local_dual_update_asyn}. These generalized algorithm updates are used to simplify the notation in the subsequent analysis.
\section*{Appendix A: Proof of Lemma \ref{lemma1}} \label{lemma2_proof}
Consider the squared 2-norm of the difference between the iterate $\x^{i}_{t}$ at time $t+1$ and an arbitrary feasible point $\bbx^i\in \ccalX^N$ and use \eqref{eq:local_primal_update_asyn_gen} to express $\bbx^{i}_{t+1}$ in terms of  $\bbx^{i}_{t}$,
\begin{align} \label{eq:primal_dist} 
\!\!\!\!\!\|\bbx^{i}_{t+1}\!\!-\!\!\bbx^i \|^2 
\!=\! \|\small\ccalP_{\ccalX}\normalfont[\bbx^{i}_{t} \!-\! \eps \nabla_{\x^i}\hat{\mathcal{L}}^i_{[\t]}\left(\!\!\{\bbx^{j}_{[t]_j},\bblam^j_t\}_{j\in n'_i} \!\right)\!]\!-\!\bbx^i \|^2.
\end{align}
\blue{where, we have utilized the compact notation defined in \eqref{d1} to substitute \eqref{eq:local_primal_update_asyn_gen} into \eqref{eq:primal_dist}. }Since $\bbx^i\in \ccalX$, utilizing non-expansive property of the projection operator in \eqref{eq:primal_dist} and expanding the square
\small
\begin{align} \label{eq:primal_sq_expand} 
\|\bbx^{i}_{t+1}-\bbx^i \|^2 \!
\leq & \|\bbx^{i}_{t}-\eps\nabla_{\bbx^i}\hat{\mathcal{L}}^i_{[\t]}\left(\{\bbx^{j}_{[t]_j},\bblam^j_t\}_{j\in n'_i}\right)-\bbx^i\|^2  \nonumber
\\ 
=&   \|\bbx^{i}_{t} - \bbx^i\|^2\nonumber
\\ 
&- 2\eps\ip{\nabla_{\bbx}\hat{\mathcal{L}}^i_{[\t]}\left(\{\bbx^{j}_{[t]_j},\bblam^j_t\}_{j\in n'_i}\right), (\bbx^{i}_{t}-\bbx^i)} \nonumber
\\
&+ \eps^2\|\nabla_{\bbx}\hat{\mathcal{L}}^i_{[\t]}\left(\{\bbx^{j}_{[t]_j},\bblam^j_t\}_{j\in n'_i}\right)\|^2.
\end{align}\normalsize
We reorder terms of the above expression such that the gradient inner product is on the left-hand side and then take summation over $i\in\ccalV$, yielding
%
%
Take the summation over nodes $i\in\ccalV$ on both sides, we get
\begin{align} \label{eq:primal_grad_neq222}
&\sum\limits_{i=1}^{N}\ip{\nabla_{\bbx}\hat{\mathcal{L}}^i_{[\t]}\left(\{\bbx^{j}_{[t]_j},\bblam^j_t\}_{j\in n'_i}\right), (\bbx^{i}_{t}-\bbx^i)}  \nonumber
\\
&\ \ \ \  \leq   \frac{1}{2\eps}\sum\limits_{i=1}^{N}\left(\|\bbx^{i}_{t}-\bbx^i\|^2 \!-\! \|\bbx^{i}_{t+1}-\bbx^i\|^2\right) \nonumber\\
&\hspace{8mm}+ \frac{\eps}{2}\sum\limits_{i=1}^{N}\norm{\nabla_{\bbx^i}\hat{\mathcal{L}}^i_{[\t]}\left(\{\bbx^{j}_{[t]_j},\bblam^j_t\}_{j\in n'_i}\right)}^2.
\end{align}

Let us use the following notation, 
\begin{align}\label{not}
\x_{[\t]}:=\left[
\x^{1}_{[\t]_1} ; \cdots
\x^{N}_{[\t]_N} 
\right],\ \text{and}\ \ 
\bblam_{t}:=\left[
\bblam^1_{t} ;\cdots 
;\bblam^N_{t}
\right].
\end{align}
Utilizing this notation, we can write
\begin{align} \label{eq:reorder}
&\ip{\nabla_{\bbx}\hat{\mathcal{L}}_{[\t]}(\x_{[\t]},\bblam_{t}), (\bbx_{t}-\bbx)}  \nonumber
\\
& \leq \!\!  \frac{1}{2\eps}\!\left(\|\bbx_{t}\!-\!\bbx\|^2 \!\!-\!\! \|\bbx_{t+1}\!\!-\!\!\bbx\|^2\right) 
\!\!+\!\! \frac{\eps}{2}\|\nabla_{\bbx}\hat{\mathcal{L}}_{[\t]}(\x_{[\t]},\bblam_{t})\|^2.
\end{align}Add and subtract $\ip{\nabla_{\bbx}\hat{\mathcal{L}}_{[\t]}(\x_{[\t]},\bblam_{t}), \bbx_{[\t]}}$ to left hand side of \eqref{eq:reorder} to obtain
\begin{align}
&\ip{\nabla_{\bbx}\hat{\mathcal{L}}_{[\t]}(\x_{[\t]},\bblam_{t}), (\bbx_{[\t]}-\bbx)}  \nonumber
\\
&\   \leq   \frac{1}{2\eps}\left(\|\bbx_{t}-\bbx\|^2 \!-\! \|\bbx_{t+1}-\bbx\|^2\right) 
+ \frac{\eps}{2}\norm{\nabla_{\bbx}\hat{\mathcal{L}}_{[\t]}(\x_{[\t]},\bblam_{t})}^2\nonumber
\\
&\ \ \ \ +\ip{\nabla_{\bbx}\hat{\mathcal{L}}_{[\t]}(\x_{[\t]},\bblam_{t}), (\bbx_{[\t]}-\bbx_{t})}. \label{eq:primal_grad_neq222222}
\end{align}
Observe now that since the functions  $f^{i}(\bbx^{i},\bbtheta^{i}_{t})$ and $\h^{i}(\{\bbx^j,\bbtheta^{j}_t\}_{j\in n'_i} )$ are convex with respect to optimization variables for any given realization of the associated random variables, therefore the stochastic Lagrangian is a convex function of $\bbx^i$ and $\bbx^j$ [cf. \eqref{eq:lagrangian}]. Hence, from the first order convex inequality and the definition of Lagrangian in \eqref{eq:stoch_lagrangian}, it holds that
\begin{align} \label{eq:primal_cvx2}
\!\!\!\!\!\hat{\mathcal{L}}_{[\t]}(\bbx_{[\t]},\bblam_{t}) \!- \!\hat{\mathcal{L}}_{[\t]}(\bbx,\bblam_{t})
\!\leq\! \ip{\nabla_{\bbx}\hat{\mathcal{L}}_{[\t]}(\bbx_{[\t]},\bblam_{t}),\! (\bbx_{[\t]}\!-\!\bbx)}.
\end{align}
Substituting the upper bound in \eqref{eq:primal_grad_neq222222} into the right hand side of  \eqref{eq:primal_cvx2} yields
\begin{align}
\hat{\mathcal{L}}_{[\t]}&(\bbx_{[\t]},\bblam_{t}) - \hat{\mathcal{L}}_{[\t]}(\bbx,\bblam_{t})    \nonumber
\\
& \leq   \frac{1}{2\eps}\left(\|\bbx_{t}-\bbx\|^2 \!-\! \|\bbx_{t+1}-\bbx\|^2\right) 
+ \frac{\eps}{2}\norm{\nabla_{\bbx}\hat{\mathcal{L}}_{[\t]}(\x_{[\t]},\bblam_{t})}^2\nonumber
\\
&\ \ \ +\ip{\nabla_{\bbx}\hat{\mathcal{L}}_{[\t]}(\x_{[\t]},\bblam_{t}), (\bbx_{[\t]}-\bbx_{t})}.\label{last2}
\end{align}
We set this analysis aside and proceed to repeat the steps in \eqref{eq:primal_dist}-\eqref{last2} for the distance between the iterate $\bblam_{t+1}^i$ at time $t+1$ and an arbitrary multiplier $\bblam^i$.
\begin{align}
\!\!\!\!\|\bblam_{t+1}^i\!\!-\!\!\bblam^i \|^2  
\!\!=\!\!\|  [ \bblam_t^i\! +\! \eps \nabla_{\bblam^i}\hat{\mathcal{L}}^i_{[\t]}\!\!\left(\!\!\{\bbx^{j}_{[t]_j},\bblam^j_t\}_{j\in n'_i}\right) ]_+\!\!-\!\!\bblam^i \|^2\label{eq:dual_dist}
\end{align}
where we have substituted \eqref{eq:local_dual_update_asyn_gen} to express $\bblam_{t+1}^i$ in terms of $\bblam_t^i$.
Using the non-expansive property of the projection operator in \eqref{eq:dual_dist} and expanding the square, we obtain
\begin{align} \label{eq:dual_sq_expand}
\| \bblam_{t+1}^i-\bblam^i\|^2 
&\leq \|\bblam_t^i + \eps \nabla_{\bblam^i}\hat{\mathcal{L}}^i_{[\t]}\left(\{\bbx^{j}_{[t]_j},\bblam^j_t\}_{j\in n'_i}\right) -  \bblam^i \|^2\nonumber
\\
& = \|\bblam_t^i - \bblam\|^2 \nonumber
\\
&\ \ +2\eps \ip{\nabla_{\bblam^i}\hat{\mathcal{L}}^i_{[\t]}\left(\{\bbx^{j}_{[t]_j},\bblam^j_t\}_{j\in n'_i}\right), (\bblam_t^i-\bblam^i)} \nonumber
\\
&\quad+ \eps^2 \|\nabla_{\bblam^i}\hat{\mathcal{L}}^i_{[\t]}\left(\{\bbx^{j}_{[t]_j},\bblam^j_t\}_{j\in n'_i}\right)\|^2 . 
\end{align}
Reorder terms in the above expression such that the gradient inner product term is on the left-hand side as
\begin{align}
\ip{\nabla_{\bblam^i}\hat{\mathcal{L}}^i_{[\t]}\left(\{\bbx^{j}_{[t]_j},\bblam^j_t\}_{j\in n'_i}\right), (\bblam_t^i-\bblam^i)} \nonumber
\\ 
&\hspace{-4.5cm}\geq \frac{1}{2\eps} \left( \| \bblam_{t+1}^i-\bblam^i\|^2 - \| \bblam_{t}^i-\bblam^i\|^2 \right) 
\!\nonumber
\\
&\hspace{-4cm}-\!\frac{\eps}{2}\|\nabla_{\bblam^i}\hat{\mathcal{L}}^i_{[\t]}\left(\{\bbx^{j}_{[t]_j},\bblam^j_t\}_{j\in n'_i}\right)\|^2.  \label{eq:dual_grad_neq}
\end{align}
Take the summation over nodes $i\in\ccalV$ so that we may write
\begin{align}
&\sum\limits_{i=1}^{N}\ip{\nabla_{\bblam^i}\hat{\mathcal{L}}^i_{[\t]}\left(\{\bbx^{j}_{[t]_j},\bblam^j_t\}_{j\in n'_i}\right), (\bblam_t^i-\bblam^i)} \nonumber
\\ 
&\hspace{1.5cm}\geq \frac{1}{2\eps} \sum\limits_{i=1}^{N}\left( \| \bblam_{t+1}^i-\bblam^i\|^2 - \| \bblam_{t}^i-\bblam^i\|^2 \right) 
\nonumber
\\
&\hspace{2cm}-\!\sum\limits_{i=1}^{N}\frac{\eps}{2}\|\nabla_{\bblam^i}\hat{\mathcal{L}}^i_{[\t]}\left(\{\bbx^{j}_{[t]_j},\bblam^j_t\}_{j\in n'_i}\right)\|^2. \label{eq:dual_grad_neq2}
\end{align}
Utilizing the notation defined in \eqref{not}, we write \eqref{eq:dual_grad_neq2} as follows
\begin{align}
&\ip{\nabla_{\bblam}\hat{\mathcal{L}}_{[\t]}(\x_{[\t]},\bblam_{t}), (\bblam_t-\bblam) }\nonumber
\\ 
&\!\!\geq \!\!\frac{1}{2\eps}\! \left( \| \bblam_{t+1}\!\!-\!\!\bblam\|^2 \!\!-\!\! \| \bblam_{t}\!\!-\!\!\bblam\|^2 \right)\! -\!\frac{\eps}{2}\|\nabla_{\bblam}\hat{\mathcal{L}}_{[\t]}(\x_{[\t]},\bblam_{t})\|^2.  \label{eq:dual_grad_neq22}
\end{align}
Note that the online Lagrangian [cf. \eqref{eq:lagrangian}] is a concave function of its Lagrange multipliers, which implies that instantaneous Lagrangian differences for fixed $\x_{[\t]}$ satisfy 
\begin{align}  \label{eq:dual_conc}
\!\!\!\!\!\!\hat{\mathcal{L}}_{[\t]}(\x_{[\t]},\bblam_{t}) \!- \!\hat{\mathcal{L}}_{[\t]}(\x_{[\t]},\bblam)
\geq \nabla_{\bblam}\hat{\mathcal{L}}_{[\t]}(\x_{[\t]},\bblam_{t})^T(\bblam_t\!-\!\bblam).
\end{align} 
By using the lower bound stated in  \eqref{eq:dual_grad_neq22} for the right hand side of \eqref{eq:dual_conc}, we can write
\begin{align} 
\hat{\mathcal{L}}_{[\t]}(\x_{[\t]},\bblam_{t})\! -\! \hat{\mathcal{L}}_{[\t]}(\x_{[\t]},\bblam)\  
&\!\!\geq\!\! \frac{1}{2\eps} \!\!\left( \| \bblam_{t+1}\!\!-\!\!\bblam\|^2 \!\!-\! \| \bblam_{t}-\bblam\|^2 \right) 
\nonumber
\\  
&\  -\!\frac{\eps}{2}\|\nabla_{\bblam}\hat{\mathcal{L}}_{[\t]}(\x_{[\t]},\bblam_{t})\|^2.\label{eq:lagrange_dual_neq}
\end{align}
We now turn to establishing a telescopic property of the instantaneous Lagrangian by combining the expressions in \eqref{last2} and \eqref{eq:lagrange_dual_neq}. To do so observe that the term $\hat{\mathcal{L}}_{[\t]}(\x_{[\t]},\bblam_{t})$ appears in both inequalities. Thus, subtracting the inequality \eqref{eq:lagrange_dual_neq} from those in \eqref{last2} followed by reordering terms yields the required result in \eqref{eq:lemma12}.
%
$\qed$
\vspace{-5mm}
\section*{Appendix B: Proof of Lemma \ref{lemma12}} \label{lemma2_proof2}

We first consider the expression in \eqref{eq:lemma12}, and expand the left-hand side using the definition of the augmented Lagrangian in \eqref{eq:stoch_lagrangian}. Doing so yields the following expression,
\vspace{-4mm}
\begin{align} 
& \sum\limits_{i=1}^{N}[f^i(\bbx^{i}_{[t]_i},\bbtheta^{i}_{[t]_i}) 
- f^i(\bbx^{i},\bbtheta^{i}_{[t]_i})] 
+\frac{\delta \eps}{2}( \norm{\bblam_{t}}^2\! -\! \norm{\bblam}^2) \nonumber
\\
&+ \sum\limits_{i=1}^{N}\left[
\ip{\bblam^{i},  \h^{i}\left(\{\bbx^{j}_{[t]_j},\bblam^j_t\}_{j\in n'_i}\right)}
\!-\! \ip{\bblam^{i}_t,  \h^{i}\left(\{\bbx^{j},\bblam^j_t\}_{j\in n'_i}\right)} \right] \nonumber
\\
&\quad \leq \frac{1}{2\eps} \big( \|\bbx_{t}\!-\!\bbx\|^2 \!- \!\|\bbx_{t+1}\!-\!\bbx\|^2 
+\|\bblam_{t}-\bblam\|^2 \!-\! \|\bblam_{t+1}\!-\!\bblam\|^2 \big)  \nonumber
\\
&\qquad	 +\frac{\eps}{2} \left(\norm{\nabla_{\bblam}\hat{\mathcal{L}}_{[\t]}(\x_{[\t]},\bblam_{t})}^2 +\| \nabla_{\bbx}\hat{\mathcal{L}}_{[\t]}(\x_{[\t]},\bblam_{t})\|^2\right)\nonumber
\\
&\quad\quad+\ip{\nabla_{\bbx}\hat{\mathcal{L}}_{[\t]}(\x_{[\t]},\bblam_{t}), (\bbx_{[\t]}-\bbx_{t})}.\label{eq:stoch_lagrangian_expand}
\end{align}
Let $\ccalF_{t}$ denotes the sigma field collecting the algorithm history which collects the information for all random quantities as  $\{\bbtheta_u, \bbx_u,\bblambda_u \}_{u < t}$.
{Note that this notation is slightly different from the standard notation in literature and  the maximum value $u$ can take is $t-1$ which is for the synchronous case}.  Note that the conditional expectation of the following term for given sigma algebra $\ccalF_{[\t]}$ is equal to  
\vspace{-5mm}\begin{align}\label{inter}
\!\!\!\!\mathbb{E}\Big[\small\!\sum\limits_{i=1}^{N}[f^i(\bbx^{i}_{[t]_i},\bbtheta^{i}_{[t]_i}) 
\!-\! f^i(\bbx^{i},\bbtheta^{i}_{[t]_i})] ~|~ \ccalF_{[\t]}\normalfont\Big]\!=\! F(\bbx_{[\t]}) 
\!-\! F(\bbx)
\end{align}
Taking the total expectation of \eqref{eq:stoch_lagrangian_expand} and utilizing the simplified expression in \eqref{inter}, we get
\begin{align} 
& \mathbb{E}[F(\bbx_{[\t]}) 
- F(\bbx)] 
+\frac{\delta \eps}{2}\mathbb{E}( \norm{\bblam_{t}}^2\! -\! \norm{\bblam}^2) \nonumber
\\
&+\!\!\sum\limits_{i=1}^{N}\mathbb{E}\left[\!\ip{
	\bblam^{i}, \h^{i}\!\left(\!\!\{\bbx^{j}_{[t]_j},\bbtheta^j_{[t]_j}\}_{j\in n'_i}  \right)\!\!}
\!-\! \ip{\bblam^{i}_t, \h^{i}\!\!\left(\{\bbx^{j},\bbtheta^j_{[t]_j}\}_{j\in n'_i} \right)\!}\!  \right] \nonumber
\\
&\leq \frac{1}{2\eps} \mathbb{E}\big( \|\bbx_{t}\!-\!\bbx\|^2 \!- \!\|\bbx_{t+1}\!-\!\bbx\|^2 
+\|\bblam_{t}-\bblam\|^2 \!-\! \|\bblam_{t+1}\!-\!\bblam\|^2 \big)  \nonumber
\\
&\qquad	 +\frac{\eps}{2} \left(\mathbb{E}\norm{\nabla_{\bblam}\hat{\mathcal{L}}_{[\t]}(\x_{[\t]},\bblam_{t})}^2 +\mathbb{E}\| \nabla_{\bbx}\hat{\mathcal{L}}_{[\t]}(\x_{[\t]},\bblam_{t})\|^2\right)\nonumber
\end{align}
\begin{align}
&\qquad+\Ex{\ip{\nabla_{\bbx}\hat{\mathcal{L}}_{[\t]}(\x_{[\t]},\bblam_{t}), (\bbx_{[\t]}-\bbx_{t})}}.\label{eq:stoch_lagrangian_expand_exp}
\end{align}	
Let us develop the upper bound on the term $\Ex{\ip{\bblam^{i}_t, \h^{i}\left(\{\bbx^{j},\bbtheta^j_{[t]_j}\}_{j\in n'_i} \right)}}$. Note that 
\begin{align}\label{first_law}
&\!\!\!\!\!\Ex{\! \ip{\bblam^{i}_t,\! \h^{i}\!\!\left(\!\{\bbx^{j},\!\bbtheta^j_{[t]_j}\!\!\}_{j\in n'_i}\! \right)\!}\!}
\!=\!\Ex{\!\Ex{\!\ip{\!\bblam^{i}_t,\! \h^{i}\!\!\left(\!\{\!\bbx^{j},\!\bbtheta^j_{[t]_j}\!\}_{j\in n'_i} \!\right)\!}\!\!~|~\!\!\ccalF_{[\t]}\!}\!}\nonumber\\
&\hspace{0.5cm}=\Ex{\ip{\bblam^{i}_t, \Ex{\left(\!\h^{i}\left(\{\bbx^{j},\bbtheta^j_{[t]_j}\}_{j\in n'_i}\! \right) \right)~|~\ccalF_{[\t]}}\!}\!} \!\leq\! 0, 
\end{align}
%
where the first equality in \eqref{first_law} holds from the law of iterated averages. The second equality holds since $\bblam^{i}_t$ is deterministic for given $\ccalF_{[\t]}$, and the third is due to the fact that for any feasible $\{\x^j\}_{j\in n'_i}$, $\mathbb{E}[\h^{i}\left(\{\bbx^{j},\bbtheta^j_{[t]_j}\}_{j\in n'_i} \right)]\leq 0$ due to the Slater's conditions (Assumption \ref{A2}), where here we use the generalized constraint definition given in \eqref{eq:const_gen} which subsumes \eqref{constraint_gen}.
Now, let's use \eqref{first_law} in the left hand side of \eqref{eq:stoch_lagrangian_expand_exp} to obtain
\begin{align} \label{eq:stoch_lagrangian_expand_exp2}
&\mathbb{E}[F(\bbx_{[\t]}) 
- F(\bbx)] 
+\frac{\delta \eps}{2}\mathbb{E}( \norm{\bblam_{t}}^2\! -\! \norm{\bblam}^2)
\\
&+\sum\limits_{i=1}^{N}\mathbb{E}\left[\ip{
	\bblam^{i}, \h^{i}\left(\{\bbx^{j}_{[t]_j},\bbtheta^j_{[t]_j}\}_{j\in n'_i}  \right)}\right] \nonumber
\\
&\leq \frac{1}{2\eps} \mathbb{E}\big( \|\bbx_{t}\!-\!\bbx\|^2 \!- \!\|\bbx_{t+1}\!-\!\bbx\|^2 
+\|\bblam_{t}-\bblam\|^2 \!-\! \|\bblam_{t+1}\!-\!\bblam\|^2 \big)  \nonumber
\\
&\ \ +\frac{\eps}{2} \left(\mathbb{E}\norm{\nabla_{\bblam}\hat{\mathcal{L}}_{[\t]}(\x_{[\t]},\bblam_{t})}^2 +\mathbb{E}\| \nabla_{\bbx}\hat{\mathcal{L}}_{[\t]}(\x_{[\t]},\bblam_{t})\|^2\right)+I,\nonumber
\end{align}
where in \eqref{eq:stoch_lagrangian_expand_exp2} we have defined $I$ on the right-hand side as
\begin{align}\label{eq:asynchrony_error}
I:= \ip{\nabla_{\bbx}\hat{\mathcal{L}}_{[\t]}(\x_{[\t]},\bblam_{t}), (\bbx_{[\t]}-\bbx_{t})}.
\end{align}
Observe that \eqref{eq:asynchrony_error} is the directional error of the stochastic gradient caused by asynchronous updates. To proceed further, an upper bound on the term $I$ is required, which we derive next. Using Cauchy Schwartz inequality, we obtain
\begin{align}\label{first}
{I}
&\leq\norm{\nabla_{\bbx}\hat{\mathcal{L}}_{[\t]}(\x_{[\t]},\bblam_{t})}\norm{(\x_{[\t]}-\x_t)}.
\end{align}
Since the maximum delay is $\tau$, we can write the difference as $\norm{\x_{t}-\x_{[\t]}}$ as sum of $\tau$ intermediate difference using triangular inequality as follows
\begin{align}\label{running_triangle_2}
\!\!\!\!\norm{\x_{t}\!-\!\x_{[\t]}}&\!\leq\!\!\!\! \sum\limits_{s=t-\tau}^{t-1}\!\!\!\norm{\x_{s+1}\!-\!\x_{{s}}}\!\leq\! \ep\!\!\!\!\sum\limits_{s=t-\tau}^{t-1}\!\!\!\!\norm{\nabla_\x {\ccalL}_{[s]}(\x_{[s]},\bblam_s)}.
\end{align}
where the inequality in \eqref{running_triangle_2} follows from the primal update in \eqref{eq:local_primal_update_asyn_gen}. For brevity, let us denote 
%
%
$B_t:=\norm{\nabla_\x {\ccalL}_{[\t]}(\x_{[\t]},\bblam_t)}$.
%
%
Note that in the definition of $B_t$, only $t$ index is  emphasized since the dual variable involved is $\bblam_t$. Substituting  upper bound obtained in \eqref{running_triangle_2} into \eqref{first}, simplyfying the notation using  the definition of $B_t$ and then taking expectation, we get 
\begin{align}\label{I_2}
\Ex{I}&\leq\ep\sum\limits_{s=t-\tau}^{t-1}\mathbb{E}\left[B_tB_s\right]\leq \frac{\ep}{2}\sum\limits_{s=t-\tau}^{t-1}{\Ex{B_t^2+B_s^2}}\nonumber
\\
&\hspace{2.4cm}= \frac{\ep}{2}\Big[\col{\tau}\cdot\Ex{B_t^2}\!\!+\!\!\!\!\sum\limits_{s=t-\tau}^{t-1}\Ex{B_s^2}\Big].
\end{align}
The second inequality in \eqref{I_2} follows directly using $ab\leq\frac{a^2+b^2}{2}$. The last equality of \eqref{I_2} is obtained  by expanding the summation. 
Next, applying the gradient norm square upper bound of \eqref{gradient_1} and \eqref{gradient_2} into \eqref{I_2}, we obtain 
\begin{align}\label{I_4}
\Ex{I}\leq &\frac{\ep}{2}\Big[\left({2}\tau{(N+M^2)} L^2 (1+\Ex{\norm{\boldsymbol{\lambda}_t}^2})\right)\nonumber
\\
&+\sum\limits_{s=t-\tau}^{t-1}\left({2(N+M^2)} L^2 (1+\mathbb{E}[\norm{\bblam_s}^2])\right)\Big].
\end{align}
Rearranging and collecting the like terms, we get
\begin{align}\label{I_5}
\Ex{I}\leq &\ep\Big[{2}\tau{(N+M^2)} L^2+\tau{(N+M^2)}L^2\Ex{\norm{\boldsymbol{\lambda}_t}^2}\!\!\nonumber
\\
&+\!\!
{(N+M^2)}L^2\sum\limits_{s=t-\tau}^{t-1}\! \mathbb{E}[\norm{\bblam_s}^2]\Big].
\end{align}
Multiplying the last term of \eqref{I_5} by $\tau$, we get
\begin{align}\label{mul_tau}
\Ex{I}\leq &\ep\Big[{2}\tau{(N+M^2)} L^2+\tau{(N+M^2)}L^2\mathbb{E}[\norm{\boldsymbol{\lambda}_t}^2]\!\!\nonumber
\\
&+\tau{(N+M^2)}L^2\sum\limits_{s=t-\tau}^{t-1} \mathbb{E}[\norm{\bblam_s}^2]\Big]
\end{align}
Let us define $K_1\!\!:=\!\!{(N\!\!+\!\!M^
	2)L^2}$, expression in \eqref{mul_tau} can be expressed as
%
\begin{align}\label{bou}
\Ex{I}\leq \ep\tau K_1\Big[2+\sum\limits_{s=t-\tau}^{t}\mathbb{E}[\norm{\bblam_s}^2]\Big].
\end{align}
Substitute the bounds developed in \eqref{gradient_1}, \eqref{gradient_2} and \eqref{bou} back into \eqref{eq:stoch_lagrangian_expand_exp2}, we get 
\begin{align}\label{eq:big_mess}
& \mathbb{E}[F(\bbx_{[\t]}) - F(\bbx)] 
+\frac{\delta \eps}{2}\mathbb{E}( \norm{\bblam_{t}}^2\! -\! \norm{\bblam}^2) \nonumber
\\
& + \sum_{i=1}^{N}\mathbb{E}\left[
\ip{\bblam^{i},  \h^{i}\left(\{\bbx^{j}_{[t]_j},\bbtheta^j_{[t]_j}\}_{j\in n'_i}  \right)}\right] \nonumber
\\
&\leq \frac{1}{2\eps} \mathbb{E}\big( \|\bbx_{t}\!\!-\!\!\bbx\|^2 \!- \!\|\bbx_{t+1}\!-\!\bbx\|^2 
\!+\!\|\bblam_{t}\!\!-\!\!\bblam\|^2 \!-\! \|\bblam_{t+1}\!-\!\bblam\|^2 \big)  \nonumber
\\
&\quad \!\!+\!\!\frac{\eps}{2}\!\! \left(\!2M\sigma_{\boldsymbol{\lambda}}^2\!\!+\!  2\delta^2\! \eps^2 \Ex{\!\| \bblambda_t \|^2\!}  \!\!+{\!\!2(N\!\!+\!\!M^2)} L^2 (1+\mathbb{E}[\norm{\bblam_t}^2])\right)\nonumber
\\
&\quad+\ep\col{\tau}K_1\Big[2+\sum\limits_{s=t-\tau}^{t}\mathbb{E}[\norm{\bblam_s}^2]\Big].
\end{align}
Collecting the like terms together and defining $K_2:=M\sigma_{\boldsymbol{\lambda}}^2+{(N\!\!+\!\!M^2)} L^2\!+\!\tau K_1$ and $K_3\!:=\!\delta^2\eps^2\!+\!{(N\!\!+\!\!M^2)}L^2$, we get
\begin{align}\label{last1}
& \mathbb{E}[F(\bbx_{[\t]}) - F(\bbx)] 
+\frac{\delta \eps}{2}\mathbb{E}( \norm{\bblam_{t}}^2\! -\! \norm{\bblam}^2) \nonumber
\\
&\quad\quad + \sum_{i=1}^{N}\mathbb{E}\left[
\ip{\bblam^{i},  \h^{i}\left(\{\bbx^{j}_{[t]_j},\bbtheta^j_{[t]_j}\}_{j\in n'_i}  \right)}\right] \nonumber
\end{align}
\begin{align}
& \leq\frac{1}{2\eps} \mathbb{E}\big( \|\bbx_{t}\!-\!\bbx\|^2 \!- \!\|\bbx_{t+1}\!-\!\bbx\|^2 
+\|\bblam_{t}-\bblam\|^2 \!-\! \|\bblam_{t+1}\!-\!\bblam\|^2 \big)  \nonumber
\\
&\quad	 +\eps\Big[K_2\!+\!K_3\mathbb{E}[\norm{\bblam_t}^2]\!+\!\frac{\tau K_1}{2}\sum\limits_{s=t-\tau}^{t}\mathbb{E}[\norm{\bblam_s}^2]\Big].
\end{align}
%
Adding $\Ex{F(\x_{[\t]})-F(\x_{t})}$ to the both sides of \eqref{last1} {and apply Lipschitz continuity of the objective, } yields
\begin{align}\label{last}
& \mathbb{E}[F(\bbx_{t}) 
- F(\bbx)] 
+\frac{\delta \eps}{2}\mathbb{E}( \norm{\bblam_{t}}^2\! -\! \norm{\bblam}^2) \nonumber
\\
& + \sum_{i=1}^{N}\mathbb{E}\left[
\ip{\bblam^{i},  \h^{i}\left(\{\bbx^{j}_{[t]_j},\bbtheta^j_{[t]_j}\}_{j\in n'_i}  \right)}\right] 
\end{align}
\begin{align}
&\leq \frac{1}{2\eps} \mathbb{E}\big( \|\bbx_{t}\!-\!\bbx\|^2 \!\!- \!\|\bbx_{t+1}\!-\!\bbx\|^2 
\!\!+\!\!\|\bblam_{t}\!-\!\bblam\|^2 \!-\! \|\bblam_{t+1}\!-\!\bblam\|^2 \big)  \nonumber
\\
&\quad	 +\!\eps\Big[\!K_2\!\!+\!\!K_3\mathbb{E}[\norm{\bblam_t}^2]\!\!+\!\frac{\tau K_1}{2}\!\!\!\! \sum\limits_{s=t-\tau}^{t}\!\!\!\!\!\mathbb{E}[\norm{\bblam_s}^2]\Big]\!\!+\!\!L_f\Ex{\norm{\x_{t}\!-\!\x_{[\t]}}}.\nonumber
\end{align}
{Now, we proceed to analyze the resulting term $L_f\Ex{\norm{\x_{t}\!-\!\x_{[\t]}}}$.} Further note that using $[\Ex{X}]^2\leq\Ex{X^2}$ for any random variable $X$, we can write
%
\begin{align}\label{der_2_1}
\!\!\!\!\!\!\!\!\!\!\mathbb{E}[\norm{\x_{t}\!\!-\!\x_{[\t]}}]\!\leq \!\!\sqrt{\!\mathbb{E}\big[\!\norm{\x_{t}\!-\!\x_{[\t]}\!}^2\!\big]\!}
\leq\!\!\Big(\mathbb{E}\big[\!\big(\!\!\!\!\sum\limits_{s=t-\tau}^{t-1}\!\!\!\!\norm{\x_{k+1}\!-\!\x_{k}}\!\!\big)^2\big]\!\!\Big)^{1/2}
\end{align}
%
Inequality in \eqref{der_2_1} follows from the triangular inequality using comparable analysis to that which yields \eqref{running_triangle_2}. Further utilizing the result $\Big(\sum\limits_{i=1}^{U}a_i\Big)^2\leq U\sum\limits_{i=1}^{U}a_i^2$, we get
\begin{align}\label{der_2_2}
\mathbb{E}[\norm{\x_{t}-\x_{[\t]}}]&\leq \Big(\tau {\sum\limits_{s=t-\tau}^{t-1}\Ex{\norm{\x_{k+1}-\x_{k}}^2}}\Big)^{1/2}.
\end{align}
Now using the upper bound for single iterate different in terms of gradient as in \eqref{running_triangle_2}, we get
\begin{align}\label{der_2_4}
\mathbb{E}[\norm{\x_{t}-\x_{[\t]}}]&\!\leq \!\Big(\!\tau\epsilon^2\!\!\! \sum\limits_{s=t-\tau}^{t-1}\!\!\!\!\mathbb{E}[\norm{\nabla_\x {\ccalL}_{[s]}(\x_{[s]},\bblam_s)}^2]\Big)^{1/2}\nonumber
\\
&\leq\! \epsilon \Big(\!\!{2}\tau({N\!\!+\!\!M^2})L^2 \!\!\!{\sum\limits_{s=t-\tau}^{t-1}\!\!\![1+\mathbb{E}[\norm{\bblam_s}^2]]}\Big)^{1/2}\!\!\!.
\end{align}
Inequality in \eqref{der_2_4} holds due to application of gradient norm bounds. From the standard inequality of $\sqrt{1+Z}\leq(1+Z)$ for all $Z\geq0$, we can write   
\begin{align}\label{der_2_5}
\mathbb{E}[\norm{\x_{t}-\x_{[\t]}}]&\leq 2\eps \tau\sqrt{K_1}\!\sum\limits_{s=t-\tau}^{t-1}[1+\mathbb{E}[\norm{\bblam_s}^2]].
\end{align}
{where we pull $2\tau$ out of square root because the product is either zero or grater than $1$.} Utilizing the upper bound of \eqref{der_2_5} for the last term in right hand side of \eqref{last} {and taking  the summation over $t=1$ to $T$, we get }
%
%
\begin{align}
& \sum_{t=1}^{T}\mathbb{E}[F(\bbx_{t}) 
- F(\bbx)] 
+\sum_{t=1}^{T}\frac{\delta \eps}{2}\mathbb{E}( \norm{\bblam_{t}}^2\! -\! \norm{\bblam}^2) \nonumber
\\
&+\sum_{t=1}^{T}\sum_{i=1}^{N}\mathbb{E}\left[
\ip{\bblam^{i},  \h^{i}\left(\{\bbx^{j}_{[t]_j},\bbtheta^j_{[t]_j}\}_{j\in n'_i}  \right)}\right] \nonumber
\end{align}
\begin{align}
&\leq \frac{1}{2\eps} \mathbb{E}\big( \|\bbx_{1}\!-\!\bbx\|^2 \! 
+\|\bblam_{1}-\bblam\|^2 \!\big)\nonumber
\\
&\quad +\eps \Big[TK_2+K_3\sum\limits_{t=1}^{T}\mathbb{E}{[\norm{\bblam_t}^2]}+\frac{\tau K1}{2} \sum\limits_{t=1}^{T}\sum\limits_{s=t-\tau}^{t}\mathbb{E}[\norm{\bblam_s}^2]\Big]\nonumber
\\
&\quad+{2\eps \tau L_f\sqrt{K_1}} \sum\limits_{t=1}^{T}\sum\limits_{s=t-\tau}^{t-1}[1+\mathbb{E}[\norm{\bblam_s}^2]]\label{eq:sum_T2}
\end{align}
In \eqref{eq:sum_T2}, we exploit the telescopic property of the summand over differences in the magnitude of primal and dual iterates to a fixed primal-dual pair $(\bbx, \bblam)$ which appears as the first term on right-hand side of \eqref{eq:sum_T3}, and the fact that the resulting expression is deterministic. By assuming the dual variable is initialized as $\bblam_1=\bb0$ and then combining the like terms together,  we can upper bound the right hand side of \eqref{eq:sum_T2}, yielding
\begin{align}
& \sum_{t=1}^{T}\mathbb{E}[F(\bbx_{t}) 
- F(\bbx)] 
+\sum_{t=1}^{T}\frac{\delta \eps}{2}\mathbb{E}( \norm{\bblam_{t}}^2\! -\! \norm{\bblam}^2) \nonumber
\\
&+\sum_{t=1}^{T}\sum_{i=1}^{N}\mathbb{E}\left[
\ip{\bblam^{i},  \h^{i}\left(\{\bbx^{j}_{[t]_j},\bbtheta^j_{[t]_j}\}_{j\in n'_i}  \right)}\right] \nonumber
\\
& \leq \frac{1}{2\eps} \mathbb{E}\big( \|\bbx_{1}\!-\!\bbx\|^2 \! 
+\|\bblam\|^2 \big)  \nonumber
\\
&\quad	 +\eps T K_2 +\epsilon K_3\sum\limits_{t=1}^{T}\mathbb{E}[\norm{\bblam_t}^2]+ {2\eps \tau T L_f \sqrt{K_1}}  \nonumber
\\
&\quad+{\eps\tau(\frac{K_1}{2}+2L_f\sqrt{K_1})}\sum\limits_{t=1}^{T}\sum\limits_{s=t-\tau}^{t}[\mathbb{E}[\norm{\bblam_s}^2]].\label{eq:sum_T3} 
\end{align}	
%
Note that in \eqref{eq:sum_T3}, in order to gather terms, an extra $\norm{\bblam_t}^2$ is added to the right-hand side. We upper bound the last term on the right-hand side of \eqref{eq:sum_T3} by considering
\begin{align}\label{lam_sum_bound2}
\!\!\!\!\sum_{t=1}^{T}\!\sum\limits_{s=t-\tau}^{t}\!\!\!\norm{\boldsymbol{\lambda}_s}^2&\!\!=\!\!\norm{\boldsymbol{\lambda}_1}^2+(\norm{\boldsymbol{\lambda}_1}^2+\norm{\boldsymbol{\lambda}_2}^2)\nonumber
\\
&\ +\!(\norm{\boldsymbol{\lambda}_1}^2+\norm{\boldsymbol{\lambda}_2}^2+\norm{\boldsymbol{\lambda}_3}^2)+\cdots\nonumber
\\
&\ +\!(\norm{\boldsymbol{\lambda}_{T-\tau}}^2\!\!+\!\!\norm{\boldsymbol{\lambda}_{T-\tau+1}}^2\!\!+\!\!\cdots+\!\!\norm{\boldsymbol{\lambda}_T}^2).
\end{align} 
The relationship in \eqref{lam_sum_bound2} then simplifies to
\begin{align}\label{lam_sum_bound_22}
\sum_{t=1}^{T}\sum\limits_{s=t-\tau}^{t}\norm{\boldsymbol{\lambda}_s}^2\leq (\tau+1)\sum_{t=1}^{T}\Ex{\norm{\boldsymbol{\lambda}_t}^2}.
\end{align} 
Utilizing this on the right hand side of \eqref{eq:sum_T3}, we get
\begin{align} \label{eq:sum_T_dual_initialize2}
\mathbb{E}\Big[&\sum_{t=1}^T [F(\bbx_t) \!- \! F(\bbx)] 
\! +\!\!\!\sum_{i=1}^{N}
\ip{\bblam^{i}, \! \sum_{t=1}^T\h^{i}\!\left(\!\{\bbx^{j}_{[t]_j},\bbtheta^j_{[t]_j}\}_{j\in n'_i}  \!\right)\!\!} \nonumber \\
&- \Big(\!\frac{\delta \eps T}{2} \!+\! \frac{1}{2\eps}\Big)\! \|\bblambda\|^2 \Big]
\leq    \frac{1}{2\eps} \|\bbx_1 \!-\! \!\bbx \|^2
+  \frac{ \eps T K}{2}\nonumber
\\
&\hspace{3.2cm}+(\eps/2) (K_4\!\!-\!\!\delta)\!\!\sum_{t=1}^{T}\!\mathbb{E}[\norm{\boldsymbol{\lambda}_t}^2]	. 
\end{align}
where $K\!\!:={\!\! 2K_2 + 4\tau L_f \sqrt{K_1}}$ and $K_4\!:=\!\! {(2K_3\!+\!(\tau+1)\tau(K_1 \!+4L_f\sqrt{K_1}))}$. 
Now selecting $\delta$ such that $(K_4-\delta)\leq0$ makes the last term on the right-hand side of the preceding expression null, so that we may write
%
\begin{align}\label{eq:lemma122_gen}
\mathbb{E}\Big[&\sum_{t=1}^T [F(\bbx_t) \!- \! F(\bbx)] 
\! +\!\!\!\sum_{i=1}^{N}\Big[
\ip{\bblam^{i},  \sum_{t=1}^T\h^{i}\left(\{\bbx^{j}_{[t]_j},\bbtheta^j_{[t]_j}\}_{j\in n'_i}  \right)}\Big] \nonumber \\
&- \Big(\!\frac{\delta \eps T}{2} \!+\! \frac{1}{2\eps}\Big)\! \|\bblambda\|^2 \Big]
\leq    \frac{1}{2\eps} \|\bbx_1 \!-\! \!\bbx \|^2+  \frac{ \eps T K}{2}.
\end{align}
Observe that Lemma \ref{lemma12} is a special case of \eqref{eq:lemma122_gen} with simplified constraint functions that only allow for pair-wise coupling of the decisions of distinct nodes (see Remark \ref{remark1}). $\qed$
\section*{Appendix C: Proof of Theorem \ref{theorem1}} \label{lemma2_proof22}
At this point, we note that the left-hand side of the expression in \eqref{eq:lemma122_gen}, and hence \eqref{eq:lemma122}, consists of three terms. The first is the accumulation over time of the global sub-optimality, which is a sum of all local losses at each node as defined in \eqref{eq:dist_stoch_opt}; the second is the inner product of the an arbitrary Lagrange multiplier $\bblam$ with the time-aggregation of constraint violation; and the last depends on the magnitude of this multiplier. We may use these later terms to define an ``optimal" Lagrange multiplier to control the growth of the long-term constraint violation of the algorithm. This technique is inspired by the approach in \cite{mahdavi2012trading,jenatton2016online}. To do so, define the \emph{augmented} dual function $\tilde{g}(\bblam)$ using the later two terms on the left-hand side of \eqref{eq:sum_T_dual_initialize2}
\begin{align} 
\tilde{g}(\bblam)\!=\!\!\sum_{i=1}^{N}
\ip{\bblam^{i},  \sum_{t=1}^T\h^{i}\left(\{\bbx^{j}_{[t]_j},\bbtheta^j_{[t]_j}\}_{j\in n'_i}  \right)}-\!\! \Big(\!\frac{\delta \eps T}{2} \!+\!\! \frac{1}{2\eps}\!\Big)\! \|\bblambda\|^2 \!.\nonumber
\end{align}
Computing the gradient and solving the resulting stationary equation over the range $\reals_{+}^M$ yields
\begin{align} \label{eq:lambda_tilde}
\tilde{\bblam}^{i}= Z(\eps)\Big[\sum_{t=1}^T \h^{i}\left(\{\bbx^{j}_{[t]_j},\bbtheta^j_{[t]_j}\}_{j\in n'_i}  \right) 
\!\Big]_+
\end{align}
for all $(i,z)\in\ccalE$, where $Z(\eps):= \frac{1}{ ( T \delta \eps + 1/\eps)}$. Substituting the selection $\bblam^i=\tilde{\bblam}^i$ defined by \eqref{eq:lambda_tilde} into \eqref{eq:lemma122_gen} results in the following expression 
%
\begin{align}\label{eq:lambda_tilde_substitute_stepsize}
\!\mathbb{E}\Bigg[\!\!&\sum_{t=1}^T [F(\bbx_t) \!\!- \!\! F(\bbx)]
\!\!+\!\!Z(\eps)\!\! \sum_{i=1}^{N}\norm{\!\Big[\!\sum_{t=1}^T \!\Big(\!\h^{i}\!\!\left(\!\!\{\bbx^{j}_{[t]_j}\!,\bbtheta^j_{[t]_j}\!\}_{j\in n'_i}  \!\right)\!\!\!\Big)\! \Big]_+}^2\!\Bigg]
\nonumber \\
&\leq    
\! \frac{1}{2\eps} \|\bbx_1 \!-\! \!\bbx \|^2
+  \frac{ \eps T K }{2}\leq    
\! \frac{\sqrt{T}}{2}\!\!\left( \|\bbx_1 \!-\! \!\bbx \|^2 \! \!
+   K\right).
\end{align}
%
The second inequality in \eqref{eq:lambda_tilde_substitute_stepsize} is obtained by selecting the constant step-size $\eps\!\! =\!\! 1/\sqrt{T}$. This result
%
%
allows us to derive both the convergence of the global objective and the feasibility of the stochastic saddle point iterates.

We first consider the average objective error sequence $\mathbb{E}[F(\bbx_t) - F(\bbx^*)]$. To do so, subtract the last term on the left-hand side of \eqref{eq:lambda_tilde_substitute_stepsize} from both sides, and note that the resulting term is non-positive. This observation allows us to omit the constraint slack term in \eqref{eq:lambda_tilde_substitute_stepsize}, which taken with the selection $\bbx=\bbx^*$ [cf. \eqref{eq:coop_stoch_opt}] and pulling the expectation inside the summand, yields
\begin{align} 
\sum_{t=1}^T \mathbb{E}[(F(\bbx_t) \!- \! F(\bbx^*)) ]
\!\! \leq    
\! \frac{\sqrt{T}}{2}\!\!\left( \|\bbx_1 \!-\! \!\bbx^* \|^2 \! \!
+   K \right) = \ccalO(\sqrt{T}),\nonumber
\end{align}
%
%
%
which is as stated in \eqref{thm1}. 
Now we turn to establishing a sublinear growth of the constraint violation in $T$, using the expression in \eqref{eq:lambda_tilde_substitute_stepsize}. Note that from the Lipschitz continuity of the objective function, we have $| F(\x_{t}) \!- \! F(\x^\star)|\leq L_f\norm{\x_{t}\!-\!\x^\star}$. An immediate consequence of this inequality is that $ F(\x_t)\! -\!  F(\x^\star)\!\geq\!\! -2L_fR$, using this in  \eqref{eq:lambda_tilde_substitute_stepsize} yields 
\begin{align} \label{eq:feasibility1}
&\mathbb{E}\Bigg[\frac{1}{ \sqrt{T} (  \delta + 1) } \sum_{i=1}^{N}\norm{\Big[\!\sum_{t=1}^T\! \Big(\h^{i}\left(\{\bbx^{j}_{[t]_j},\bbtheta^j_{[t]_j}\}_{j\in n'_i}  \right)\Big) \Big]_+}^2\!\Bigg] \nonumber
\\
&\qquad \qquad \leq    \!
\! \frac{\sqrt{T}}{2}\!\!\left( \|\bbx_1 \!-\!\! \bbx^* \|^2 \! \!\!
+ \! \!K \right) +2TL_fR . 
\end{align}
which, after multiplying both sides by $2 \sqrt{T} (  \delta + 1)$ yields
\begin{align} \label{eq:feasibility2}
&\mathbb{E}\!\Bigg[\sum_{i=1}^{N}\norm{\Big[\!\sum_{t=1}^T\! \Big(\h^{i}\left(\{\bbx^{j}_{[t]_j},\bbtheta^j_{[t]_j}\}_{j\in n'_i}  \right)\Big) \Big]_+}^2\! \Bigg]
\\
&\qquad\qquad\!\!\! \leq    \! 
\!\Big(2 \sqrt{T} (  \delta + 1) \! \Big)\! \Big(\!\frac{\sqrt{T}}{2}\!\!\left( \|\bbx_1 \!-\! \!\bbx^* \|^2 \! \!
+  \! \!K \right) 
\!\!+2TL_fR\Big) \; . \nonumber
%
\end{align}
We complete the proof by noting that the square of the network-in-aggregate constraint violation $\sum_{i=1}^{N}\norm{\Big[\!\sum_{t=1}^T\! \Big(\h^{i}\left(\{\bbx^{j}_{[t]_j},\bbtheta^j_{[t]_j}\}_{j\in n'_i}  \right)\Big) \Big]_+}^2\!$ upper bounds the square of individual proximity constraint violations since it is a sum of positive squared terms, i.e.,
%
\small
\begin{align}\label{eq:individual_proximity_lower_bd}
&\mathbb{E}\Bigg[\sum_{i=1}^{N}\norm{\Big[\!\sum_{t=1}^T\! \Big(\h^{i}\left(\{\bbx^{j}_{[t]_j},\bbtheta^j_{[t]_j}\}_{j\in n'_i}  \right)\Big) \Big]_+}^2\!\Bigg] \!\nonumber
\\
&\quad\quad \geq \! \mathbb{E}\Bigg[\Big[\!\sum_{t=1}^T\!\left( h_{ij}(\bbx^{i}_{[t]_i},\bbx^{j}_{[t]_j},\bbtheta^i_{[t]_i},\bbtheta^j_{[t]_j})\!-\!\gamma_{ij}\right) \Big]^2_+\Bigg]
\end{align}\normalsize
\colb{where we utilized the notation defined in \eqref{constraint_gen} and the inequality that the norm square of a vector is always greater than square of  each element of the vector. The inequality in \eqref{eq:individual_proximity_lower_bd} is true for any arbitrary $ij$ in the right hand side}.

Thus the right-hand side of \eqref{eq:individual_proximity_lower_bd} may be used in place of the left-hand side of \eqref{eq:feasibility2}, implying that 
\small
\begin{align} \label{eq:feasibility3}
&\mathbb{E}\Big[\Big[\!\sum_{t=1}^T\! h_{ij}\left(\bbx^{i}_{[t]_i},\bbx^{j}_{[t]_j},\bbtheta^i_{[t]_i},\bbtheta^j_{[t]_j}\right) \Big]^2_+\Bigg]
\\
&\!\!\!\qquad\qquad \leq    \! 
\!\Big(2 \sqrt{T} (  \delta + 1) \! \Big)\! \Big(\!\frac{\sqrt{T}}{2}\!\!\left( \|\bbx_1 \!-\! \!\bbx^* \|^2 \! \!
+  \! \!K \right) 
\!\!+2TL_fR\Big) \; .\nonumber
%
\end{align}\normalsize 

%

In order to present the results for the special case discussed in \eqref{eq:coop_stoch_opt}, compute the square root of both sides of \eqref{eq:feasibility3} and take the summation over all $(i,j)\in \mathcal{E}$ to conclude \eqref{eq:theorem12}.
$\qed$
\vspace{-5mm}

\bibliographystyle{IEEEtran}
\bibliography{bibliography}

   \end{document}

%% file: my_sections.tex
\usepackage{needspace}

